\documentclass[a4paper]{article}
\usepackage[english]{babel}
\usepackage{amsthm}   
\usepackage{amsmath}  
\usepackage{amssymb}  
\usepackage{mathabx}  
\usepackage{graphicx} 
\usepackage{fancyhdr} 
\usepackage{setspace} 
\usepackage{lineno}	  
\usepackage{lastpage} 
\usepackage{booktabs} 
\usepackage{rotating} 
\usepackage{multirow}	
\usepackage{caption,subcaption}  
\usepackage[sort&compress,square,comma,numbers]{natbib}   
\usepackage{todonotes}
\usepackage{custom}
\usepackage{tikz}
\usepackage{float}
\newcommand*\circled[1]{\tikz[baseline=(char.base)]{
		\node[shape=circle,draw,inner sep=1pt] (char) {#1};}}

\newcommand{\mytitle}{Optimal operating policies for organic Rankine cycles for waste heat recovery under transient conditions}
\newcommand{\myshorttitle}{Optimal operating policies for ORCs}
\newcommand{\myauthor}{Yannic Vaupel, Wolfgang R. Huster, Adel Mhamdi and Alexander Mitsos$^{*}$} 
\newcommand{\myauthorshort}{Y. Vaupel et int., A. Mitsos}
\author{\myauthor}
\usepackage[colorlinks,linkcolor=blue,citecolor=blue,urlcolor=blue,pdftitle={\mytitle},pdfauthor={You and Alexander Mitsos}]{hyperref} 

{
	\fancyhead[L]{\myshorttitle}
	\fancyhead[R]{\small{\textit{\the\day.\the\month.\the\year}}}
	\fancyfoot[L]{Submitted to \textit{Energy} on March 5\textsuperscript{th}, 2020}
	\fancyfoot[R]{Page \thepage\ of \pageref*{LastPage}}
	\cfoot{}
}

\fancypagestyle{firststyle}
{
	 \fancyhead[L]{\copyright \small{\textit{\myauthorshort}}}
	 \fancyhead[C]{}
   \fancyhead[R]{Page \thepage\ of \pageref*{LastPage}}
	 \fancyfoot[L]{\footnotesize \rule{0.25\textwidth}{0.4pt} \newline Corresponding author: $^*$A. Mitsos \newline 
	                                                          Process Systems Engineering (AVT.SVT), RWTH Aachen University \newline
                                                            Forckenbeckstra{\ss}e 51, D - 52074 Aachen \newline
																														E-mail: amitsos@alum.mit.edu}
   \fancyfoot[C]{}
	 \fancyfoot[R]{}
}

\addto\captionsenglish{%
}

\theoremstyle{remark} 

\newcommand{\SVT}{
\begin{footnotesize}
Process Systems Engineering (AVT.SVT), \\RWTH Aachen University, 52074 Aachen, Germany
\end{footnotesize}\\
}

\newcommand{\JARAENERGY}{
\begin{footnotesize}
JARA-ENERGY, 52056 Aachen, Germany
\end{footnotesize}\\
}

\newcommand{\IEK}{
\begin{footnotesize}
Institute of Energy and Climate Research: Energy Systems Engineering (IEK-10), Forschungszentrum J{\"u}lich GmbH, 52425 J{\"u}lich, Germany.
\end{footnotesize}\\
}

\usepackage{xcolor} 

\definecolor{rwth}{RGB}{0,84,159}
\definecolor{rwth-75}{RGB}{ 64 127 183}
\definecolor{rwth-50}{RGB}{142 186 229}
\definecolor{rwth-25}{RGB}{199 221 242}
\definecolor{rwth-10}{RGB}{232 241 250}

\definecolor{black}   {RGB}{  0   0   0}
\definecolor{black-75}{RGB}{100 101 103}
\definecolor{black-50}{RGB}{156 158 159}
\definecolor{black-25}{RGB}{207 209 210}
\definecolor{black-10}{RGB}{236 237 237}

\definecolor{magenta}   {RGB}{227   0 102}
\definecolor{magenta-75}{RGB}{233  96 136}
\definecolor{magenta-50}{RGB}{241 158 177}
\definecolor{magenta-25}{RGB}{249 210 218}
\definecolor{magenta-10}{RGB}{253 238 240}

\definecolor{yellow}   {RGB}{255 237   0}
\definecolor{yellow-75}{RGB}{255 240  85}
\definecolor{yellow-50}{RGB}{255 245 155}
\definecolor{yellow-25}{RGB}{255 250 209}
\definecolor{yellow-10}{RGB}{255 253 238}

\definecolor{petrol}   {RGB}{  0  97 101}
\definecolor{petrol-75}{RGB}{ 45 127 131}
\definecolor{petrol-50}{RGB}{125 164 167}
\definecolor{petrol-25}{RGB}{191 208 209}
\definecolor{petrol-10}{RGB}{230 236 236}

\definecolor{turkis}   {RGB}{  0 152 161}
\definecolor{turkis-75}{RGB}{  0 177 183}
\definecolor{turkis-50}{RGB}{137 204 207}
\definecolor{turkis-25}{RGB}{202 231 231}
\definecolor{turkis-10}{RGB}{235 246 246}

\definecolor{grun}   {RGB}{ 87 171  39}
\definecolor{grun-75}{RGB}{141 192  96}
\definecolor{grun-50}{RGB}{184 214 152}
\definecolor{grun-25}{RGB}{221 235 206}
\definecolor{grun-10}{RGB}{242 247 236}

\definecolor{maigrun}   {RGB}{189 205   0}
\definecolor{maigrun-75}{RGB}{208 217  92}
\definecolor{maigrun-50}{RGB}{224 230 154}
\definecolor{maigrun-25}{RGB}{240 243 208}
\definecolor{maigrun-10}{RGB}{249 250 237}

\definecolor{orange}   {RGB}{246 168   0}
\definecolor{orange-75}{RGB}{250 190  80}
\definecolor{orange-50}{RGB}{253 212 143}
\definecolor{orange-25}{RGB}{254 234 201}
\definecolor{orange-10}{RGB}{255 247 234}

\definecolor{rot}   {RGB}{204   7  30}
\definecolor{rot-75}{RGB}{216  92  65}
\definecolor{rot-50}{RGB}{230 150 121}
\definecolor{rot-25}{RGB}{243 205 187}
\definecolor{rot-10}{RGB}{250 235 227}

\definecolor{bordeaux}   {RGB}{161  16  53}
\definecolor{bordeaux-75}{RGB}{182  82  86}
\definecolor{bordeaux-50}{RGB}{205 139 135}
\definecolor{bordeaux-25}{RGB}{229 197 192}
\definecolor{bordeaux-10}{RGB}{245 232 229}

\definecolor{violett}   {RGB}{ 97  33  88}
\definecolor{violett-75}{RGB}{131  78 117}
\definecolor{violett-50}{RGB}{168 133 158}
\definecolor{violett-25}{RGB}{210 192 205}
\definecolor{violett-10}{RGB}{237 229 234}

\definecolor{lila}   {RGB}{122 111 172}
\definecolor{lila-75}{RGB}{155 145 193}
\definecolor{lila-50}{RGB}{188 181 215}
\definecolor{lila-25}{RGB}{222 218 235}
\definecolor{lila-10}{RGB}{242 240 247}

\newif\ifshowmytodo
\showmytodotrue

\makeatletter
\let\@addpunct\@gobble
\g@addto@macro{\thm@space@setup}{\thm@headpunct{}} 
\makeatother

\renewenvironment{abstract}{\noindent\textbf{Abstract:}}{}

\begin{document}

\thispagestyle{firststyle}
\begin{flushleft}\begin{large}\textbf{\mytitle}\end{large} \end{flushleft}
\myauthor  
\begin{flushleft}

\JARAENERGY
\SVT
\IEK

\end{flushleft}

\begin{abstract}
	Waste heat recovery for trucks via organic Rankine cycle is a promising technology to reduce fuel consumption and emissions. 
	As the vehicles are operated in street traffic, the heat source is subject to strong fluctuations.
	Consequently, such disturbances have to be considered to enable safe and efficient operation.
	Herein, we find optimal operating policies for several representative scenarios by means of dynamic optimization and discuss the implications on control strategy design.
	First, we optimize operation of a typical driving cycle with data from a test rig.
	Results indicate that operating the cycle at minimal superheat is an appropriate operating policy.
	Second, we consider a scenario where the permissible expander power is temporarily limited, which is realistic in street traffic.
	In this case, an operating policy with flexible superheat can reduce the losses associated with operation at minimal superheat by up to 53\% in the considered scenario.
	As the duration of power limitation increases, other constraints might become active which results in part of the exhaust gas being bypassed, hence reduced savings.
\end{abstract}

\section{Introduction} \label{sec:introduction}
Due to increasing fuel prices and tighter emission regulations, waste heat recovery (WHR) from vehicles with internal combustion engines (ICE) in street traffic has become an increasingly viable option \cite{Sprouse2013,Hoang2018}.
Thermal energy is typically recovered from the exhaust gas line or exhaust gas recirculation, as the exhaust gas has higher exergy than the cooling water \cite{Fu2013}, and used as a heat source for a power cycle, e.g., a bottoming organic Rankine cycle (ORC).
\\
\\
ORCs are a proven technology in a variety of applications with low- to medium-temperature heat sources, e.g., geothermal brine \cite{Ghasemi2013} or solar thermal energy \cite{Quoilin2011b} and industrial waste heat recovery \cite{Campana2013}.
Moreover, ORC technology has been proposed in other transport applications such as marine applications \cite{Singh2016} or trains powered by ICEs \cite{Peralez2015}.
An exhaustive overview of the potential applications of ORCs is beyond the scope of this manuscript and we refer the reader to \cite{Tchanche2011}.
\\
\\
Though the WHR technology has also been proposed for passenger vehicles \cite{Horst2013,Horst2014,Boretti2012}, heavy-duty trucks are more suitable as they are operated on long distances and the additional weight of the WHR system does not have a profound effect on fuel economy.
A discussion of the various ORC configurations proposed is beyond the scope of this manuscript and we refer the reader to the reviews in \cite{Sprouse2013,Hoang2018,Lion2017,Xu2019}.
Potential fuel savings ranging from 5\% up to 10\% have been reported in simulation-based studies \cite{Sprouse2013,Lion2017}.
However, a discrepancy between fuel savings estimated in simulation-based studies and fuel savings realized in experiments is noted in \cite{Xu2019}.
In most of the aforementioned applications, the ORC system is operated either with a nearly constant heat source (geothermal), a slowly varying heat source that is predictable (solar-thermal) or with long periods of steady-state operation (ships and trains).
Heavy-duty diesel trucks, however, are operated under highly transient heat source conditions due to their use in street traffic.
\\
\\
Available publications on ORC applications with non-transient conditions typically consider optimization of steady-state operating points \cite{Dai2009,Schweidtmann2019}, some including design considerations \cite{Huster2017} and working fluid (WF) selection \cite{Wei2007,Shengjun2011,Huster2019}.
Sometimes even fluid mixtures are designed \cite{Lampe2014,Lee2017,Huster2020}.
A possibility to account for mild variable operating conditions in design optimization is by clustering operating points \cite{Mondejar2017,Schilling2018,Schilling2019} and considering off-design behavior with stationary \cite{Ghasemi2013} or simplified dynamic models \cite{Manente2013,Tillmanns2019}.
The publications in \cite{Schilling2018,Schilling2019,Tillmanns2019} consider WHR for a heavy-duty truck, whereas \cite{Mondejar2017} considers a marine application and \cite{Manente2013} consider a geothermal application. 
While steady-state models are used in most of the aforementioned publications, dynamic models might be required depending on the time scales on which changes in the inputs and disturbances occur related to the system inertia.
Thus, most publications reporting dynamic ORC models consider waste heat recovery in diesel-trucks \cite{Wei2008,Quoilin2011b,Feru2013,Seitz2016b,Koppauer2017,Xu2017,Wang2017,Huster2018} but exceptions, e.g., \cite{Eller2019}, where a geothermal ORC system is considered, exist.
\\
\\
Many interactions between exhaust gas, ORC, cooling water cycle and engine exist \cite{Grelet2016} and the significance of dynamic effects on optimal system operation is widely accepted \cite{Espinosa2011,Tona2012,Horst2014,Peralez2017}.
To achieve efficient system operation, maximizing time in power production mode, i.e., maintain sufficient superheating to allow for turbine operation, is of paramount importance \cite{Xie2013}.
The effect of transient exhaust gas conditions is considered in \cite{Xu2017a} for an ORC system in a truck with two parallel heat exchangers utilizing the tailpipe exhaust gas and exhaust gas recirculation.
The authors examine three strategies for set-point generation and find that a fuzzy logic strategy with flexible superheat exhibits the best performance.
The work is extended with respect to real-time application in \cite{Xu2020}.
\\
\\
In order to adequately control the WHR system, understanding optimal system operation is crucial.
Many of the numerous contributions on control design for ORCs consider following predefined set-point trajectories which are often obtained from steady-state optimization, neglecting dynamic effects \cite{Yebi2017}.
The majority of publications on steady-state operation of ORCs finds that operation at minimal superheating is desirable \cite{Yamamoto2001}.
Consequently, this notion is adopted in many control related publications where the control strategy aims at maintaining a fixed superheat \cite{Quoilin2011b,Peralez2013,Hernandez2016,Seitz2018,Vaupel2019b} or vapor quality \cite{Feru2014b}.
However, in \cite{Ghasemi2013}, we demonstrated that optimal off-design operation of air-cooled geothermal power plants mandates varying the superheat as a function of the ambient temperature.
The literature review above indicates that, although many publications on control strategies for ORCs operated in a transient setting are available, a relevant research gap exists with respect to the understanding of economically optimal dynamic ORC operation.
\\
\\
In this contribution, we address this gap by assessing whether the notion of operation at minimal superheat is optimal for an ORC for WHR in street traffic.
Therefore, we apply dynamic optimization to the WHR system subject to transient heat source conditions. From the optimization results, we infer optimal operating policies.
The dynamic optimization problem that we solve is closely related to the nonlinear model predictive control (NMPC) formulation, which has been applied to WHR in various publications \cite{Petr2015,Peralez2015,Liu2017,Yebi2017}.
Indeed, we solve a similar optimal control problem (OCP) which considers the full length of the respective scenario and assumes full knowledge of the heat source.
Thus, it represents an upper bound on the performance of NMPC and allows to draw conclusions for control strategy design.
We consider two exemplary cases whose characteristics occur in street traffic.
First, we consider exhaust data used for model validation in \cite{Huster2018}.
It was recorded on a test rig and consists of parts of the World Harmonized Transient Cycle (WHTC).
For simplicity, we assume that no operational restrictions, beyond safety constraints, are imposed on the system.
\\
\\
Second, we examine a scenario where the permissible turbine power is temporarily limited. 
This scenario typically occurs in street traffic when the engine torque is negative or below a certain threshold.
Though the considered system includes a battery, the permissible turbine power can be limited by the operational constraints of a battery.	
The dynamic optimization problems are solved using the open-source dynamic optimization tool DyOS \cite{Caspari2019}.
\\
\\
The remainder of this manuscript is structured as follows. We provide a brief presentation of the examined system and the model in Sec.~\ref{sec:model}, followed by a presentation of the optimization procedure in Sec.~\ref{sec:optProblem}. In Sec.~\ref{sec:nominal}, we examine optimal operation for a typical driving cycle followed by a detailed examination of a scenario where the expander power is temporarily limited (Sec.~\ref{sec:activeConstraints}). We discuss the results of the case studies and the implications on control strategy design in Sec.~\ref{sec:discussion} and present our conclusions in Sec.~\ref{sec:conclusion}.

\section{Process model} \label{sec:model}
\begin{figure}[h]
\centering
\includegraphics[width=0.9\linewidth]{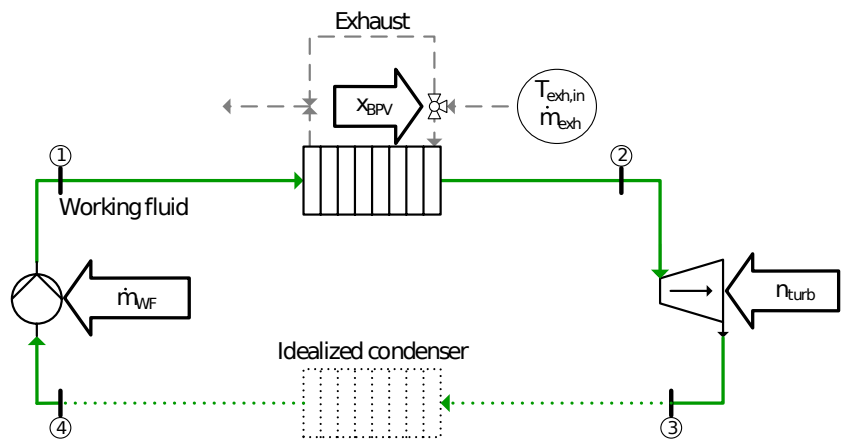}
\caption{Topology of the examined system. The WF is indicated by the solid green line and the exhaust gas by the dashed gray line. The DOFs are indicated by arrows and the disturbances by the circle. The condenser is represented with an idealized model that assumes operation at ambient pressure and fixed subcooling.}
\label{fig:topology}
\end{figure}

The system under investigation (Fig.~\ref{fig:topology}) is an ORC for WHR in a heavy-duty diesel truck operated in street traffic.
The liquid WF ethanol is compressed in a pump to a high pressure level \circled{4}$\to$\circled{1} and then evaporated and superheated in a heat exchanger \circled{1}$\to$\circled{2}.
The heat source is the exhaust gas of the diesel truck which can be bypassed through the exhaust bypass proportional control valve.
Consequently, the WF is expanded in a turbine \circled{2}$\to$\circled{3}, which is connected to an electric generator.
The WF is then condensed and subcooled in a condenser \circled{3}$\to$\circled{4}, for which a dedicated cooling cycle is employed.
\\
\\
The model representing the WHR system is based on the validated test rig model from \cite{Huster2018} and implemented in Modelica.
The evaporator is modeled using a moving boundary (MB) model and the turbine and pump are modeled using pseudo steady-state models.
The isentropic and mechanical efficiencies of the turbine depend on the fluid conditions at the inlet, the low pressure and the turbine speed.
For the pump, we assume constant isentropic and mechanic efficiencies of 0.9, respectively.
The thermodynamic properties of the WF are modeled using an implementation of the Helmholtz equation of state for ethanol \cite{Schroeder2014}.
\\
\\
We follow the common practice of focusing on the high pressure side of the system, i.e., we omit an elaborate condenser model as it does not significantly influence the high pressure part \cite{Peralez2017}.
Instead, we assume that the condenser operates at ambient pressure and the WF leaves the condenser as a subcooled liquid with fixed subcooling.
Under these assumptions, the turbine rotational speed does not influence the high pressure side in our model  and solely serves for optimizing turbine efficiency \cite{Tona2015}.
Further, we neglect heat losses in the pipes which connect the process units.
\\
\\
The desired mode of operation of the WHR system, which we refer to as ``nominal operating mode'' \cite{Vaupel2019}, describes the situation where the WF enters the evaporator as a subcooled liquid, leaves as superheated vapor and is expanded in the turbine where power is produced.
As the WF is a wet-expanding fluid, a certain level of superheat is required.
In this manuscript, the WHR system is always operating in nominal operating mode.

\section{Optimization problem}\label{sec:optProblem}
We solve dynamic optimization problems of the following type to find optimal operating policies:
\begin{align}
\mathrm{min}\quad & \Phi\left(\boldsymbol{x}\left(t_f\right)\right) \label{eq:obj} \\
\mathrm{s.\, t.}\quad & \dot{\boldsymbol{x}}\left(t\right)=\boldsymbol{f}\left(\boldsymbol{x}\left(t\right),\boldsymbol{y}\left(t\right),\boldsymbol{u}\left(t\right),\boldsymbol{d}\left(t\right)\right) \label{eq:differential}\\
& \boldsymbol{0}=\boldsymbol{g}\left(\boldsymbol{x}\left(t\right),\boldsymbol{y}\left(t\right),\boldsymbol{u}\left(t\right),\boldsymbol{d}\left(t\right)\right) \label{eq:algebraic}\\
& \boldsymbol{x}\left(t=0\right)=\boldsymbol{x}_{0} \label{eq:IC}\\
& p^{*,min}\leq p^*\left(t\right)\leq p^{*,max}\label{eq:pSafety}\\
& T_{WF,evap,out}^*\left(t\right)\leq T_{WF,evap,out}^{*,max}\label{eq:tSafety}\\
& \Delta T_{sup}^{min}\leq \Delta T_{sup}\left(t\right)\label{eq:supSafety} \\
& \dot{m}_{WF,in}^{*,min}\leq \dot{m}_{WF,in}^*\left(t\right)\leq \dot{m}_{WF,in}^{*,max} \label{eq:mInConstraints}\\
& n_{turb}^{*,min}\leq n_{turb}^*\left(t\right)\leq n_{turb}^{*,max} \label{eq:nTurbConstraint} \\
& x_{BPV}^{min}\leq x_{BPV}\left(t\right)\leq x_{BPV}^{max}\, . \label{eq:xBpvConstraint}
\end{align}
$\Phi$ is a Mayer-type objective function, i.e., it is evaluated at final time $t_f$.
The differential equations $\boldsymbol{f}$ and the algebraic equations $\boldsymbol{g}$ are specified in \eqref{eq:differential}-\eqref{eq:algebraic}, where $\boldsymbol{x}$ is the vector of differential states, $\boldsymbol{y}$ the vector of algebraic states, $t$ is the time, $\boldsymbol{u}$ are the inputs to the model determined by the optimizer and $\boldsymbol{d}$ the disturbances, i.e., the exhaust gas conditions.
The initial conditions $\boldsymbol{x}_0$ are specified in \eqref{eq:IC}.
Safety-related path constraints are specified in \eqref{eq:pSafety}-\eqref{eq:supSafety}.
Asterisks indicate quantities that are scaled in the same manner as in \cite{Huster2018} for confidentiality reasons.
$\Delta T_{sup}$ is the WF superheat where $\Delta$ indicates a temperature difference, $T_{WF,evap,out}^*$ is the WF evaporator outlet temperature and $p^*$ is the high pressure.
\\
The degrees of freedom (DOF) for optimization are the WF fluid mass flow $\dot{m}_{WF,in}^*$, the turbine rotational speed $n_{turb}^*$ and the exhaust gas valve position $x_{BPV}$, which are box-constrained in \eqref{eq:mInConstraints}-\eqref{eq:xBpvConstraint}.
\\
\\
A list of the lower and upper bounds for the path constraints and the DOF is provided in Table~\ref{tab:bounds}.
The constraint on minimal superheat protects the turbine from damage due to droplet formation.
The lower bound on pressure reflects a minimal pressure ratio.
As the condenser is assumed to operate at ambient pressure, this constraint can be directly expressed for the high pressure side.
The upper bound on $p^*$ ensures safe operation as does the maximum WF outlet temperature which prevents  WF degradation.
We do not specify a lower bound on WF outlet temperature.
However, an effective lower bound is provided at any time through the minimal superheat constraint added with the saturation temperature of the pressure at that time.
The lower and upper bounds on $\dot{m}_{WF,in}^*$ and $n_{turb}^*$ reflect the limits in which the model is valid \cite{Huster2018} and the lower bound on $x_{BPV}$ is included to avoid simulation failure due to very small exhaust gas mass flows.
\begin{table}[h!]
\centering
\begin{tabular}{lcrrlcrr}
	\toprule
	\multicolumn{4}{c}{Path} &\multicolumn{4}{c}{DOF} \\
	Variable & Unit & LB & UB & Variable & Unit & LB & UB \\
	\cmidrule(lr){1-4}
	\cmidrule(lr){5-8}
	$\Delta T_{sup}$ & K & 10 & - & $\dot{m}_{WF,in}^*$ & - & 0.0073 & 0.0363 \\
	$T_{WF,evap,out}^*$ & - & - & 0.8719 & $n_{turb}^*$ & - & 0.73 & 1.09 \\
	$p^*$& - & 0.3 & 1.5 & $x_{BPV}$ & - & 0.05 & 1.00 \\
	\bottomrule
\end{tabular}
\caption{Bounds of path constraints and DOF.}
\label{tab:bounds}
\end{table}
\\
\\
In the following sections, we compare two operating policies.
We assess the examined policies using the resulting net average power $P_{net,av}^*$, which is defined as follows
\begin{align}
P_{net,av}^*=\frac{\int_{t_{0}}^{t_{f}}\left(P_{turb}^*\left(t\right)-P_{pump}^*\left(t\right)\right)\, \mathrm{d}t}{t_f-t_0}
\end{align}
where $P_{turb}^*$ is the scaled turbine power and $P_{pump}^*$ is the scaled pump power.
\\
\\
First, we assess a policy that aims at maintaining minimal superheat while using $n_{turb}^*$ to optimize turbine efficiency.
We refer to this strategy, which is a standard approach in literature, as MSH (minimal superheat).
Fixing the superheat to a minimal value is infeasible as it often resulted in integration failure.
Furthermore no unique solution that provides minimal superheat exists, due to the availability of the exhaust gas bypass valve.
Thus, we use an optimization-based approach.
The objective can be expressed as 
\begin{align}
\Phi_{1}\left(t_f\right) = \int_{t_{0}}^{t_{f}}\left(\Delta T_{sup}\left(t\right)-\Delta T_{sup}^{min}\right)^2\, \mathrm{d}t\, .
\end{align}
As we assess the examined operating policies based on $P_{net,av}^*$, we introduce specific measures that ensure that we find the minimal superheat strategy with highest $P_{net,av}^*$ in Sec.~\ref{sec:nominal} and Sec.~\ref{sec:activeConstraints}.
\\
\\
Second, we examine the thermodynamically optimal policy, i.e., maximizing the net work without consideration of a desired superheat, which we refer to as FSH (flexible superheat).
This also corresponds to the economically optimal policy, given that all produced power can be utilized, and can be expressed as
\begin{align}
\Phi_{2}\left(t_f\right)=-\int_{t_{0}}^{t_{f}}P_{net}\left(t\right)\, \mathrm{d}t\, . \label{eq:obj41}
\end{align}
\\
\\
We assume the validated model to represent the real system behavior herein.
Since we are interested in understanding how to best operate the system, i.e., with FSH or MSH, mismatch between the system and our model is a minor concern.
However, when the model is used for the control of a physical system, considerations regarding plant-model mismatch are required.
A potential remedy could be the addition of a disturbance model to achieve offset-free model predictive control as practiced in \cite{Rathod2019}.
\\
\\
For all scenarios, the initial state of the system $\boldsymbol{x}_0$ is specified to the economically optimal steady-state, indicated as $\boldsymbol{x}_{ss}^{opt}$, for the heat source conditions at $t=0~\mathrm{s}$
\begin{align}
\boldsymbol{x}_{0} = \boldsymbol{x}_{ss}^{opt}\left({\boldsymbol{d}\left(t=0\right),\boldsymbol{u}_{ss}^{opt}}\right)\, .
\end{align}
We determine $\boldsymbol{x}_{ss}^{opt}$ in an a-priori optimization and, as expected, it corresponds to operation with minimal superheat.
\\
\\
The dynamic optimization problems are solved using direct single shooting \cite{Biegler2010} with the open-source software DyOS \cite{Caspari2019}.
The model is linked to DyOS through the functional mock-up interface (FMI).
The sensitivities are calculated through sensitivity integration with the integrator sLimex \cite{Schlegel2004} and the NLPs are solved with SNOPT \cite{Gill2005}.
All degrees of freedom are discretized on a piecewise linear continuous grid which is determined by the grid adaption algorithm described in \cite{Schlegel2005}.
\section{Optimal operation considering only safety constraints} \label{sec:nominal}
\begin{figure}[h]
\centering
\begin{subfigure}[h]{0.48\linewidth}
	\includegraphics[width=\textwidth]{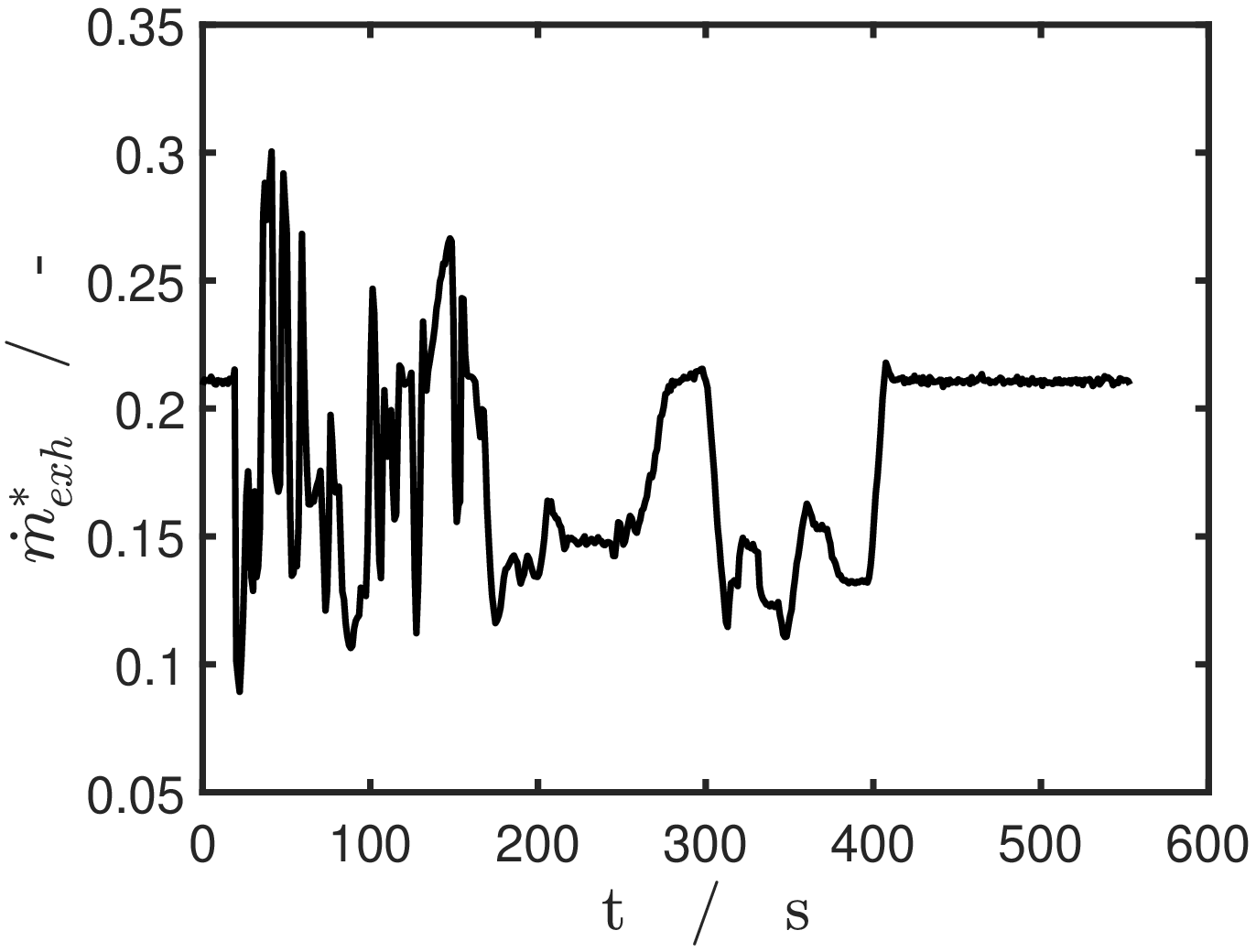}
	\caption{Exhaust gas mass flow}
	\label{fig:mExhWhtc}
\end{subfigure}
~
\begin{subfigure}[h]{0.48\linewidth}
	\includegraphics[width=\textwidth]{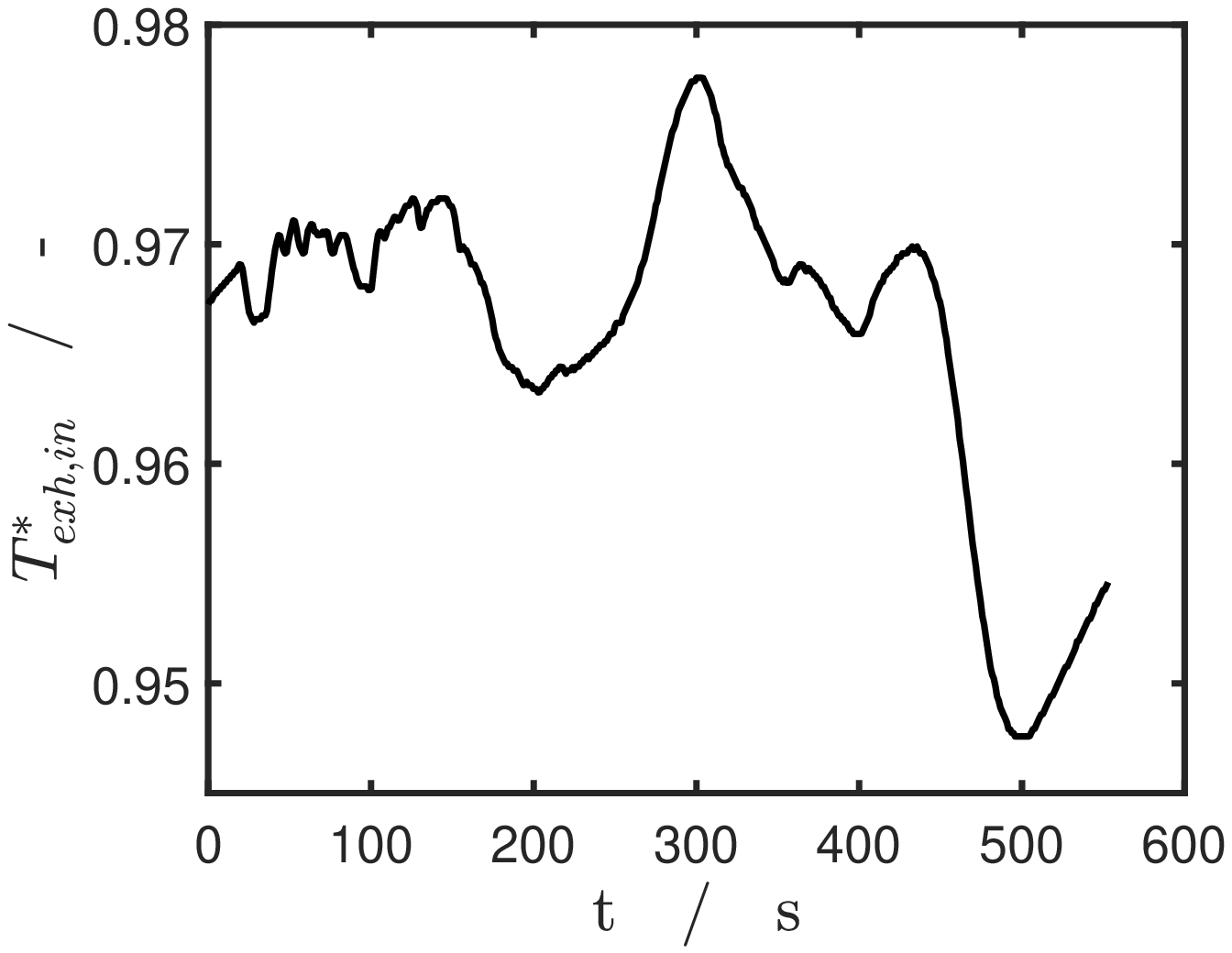}
	\caption{Exhaust gas inlet temperature}
	\label{fig:tExhWhtc}
\end{subfigure}
\caption{Heat source data for the WHTC taken from the test rig described in \cite{Huster2018}}\label{fig:heatSourceWhtc}
\end{figure}

\begin{figure}[h!]
\centering
\begin{subfigure}[h]{0.45\linewidth}
	\includegraphics[width=\textwidth]{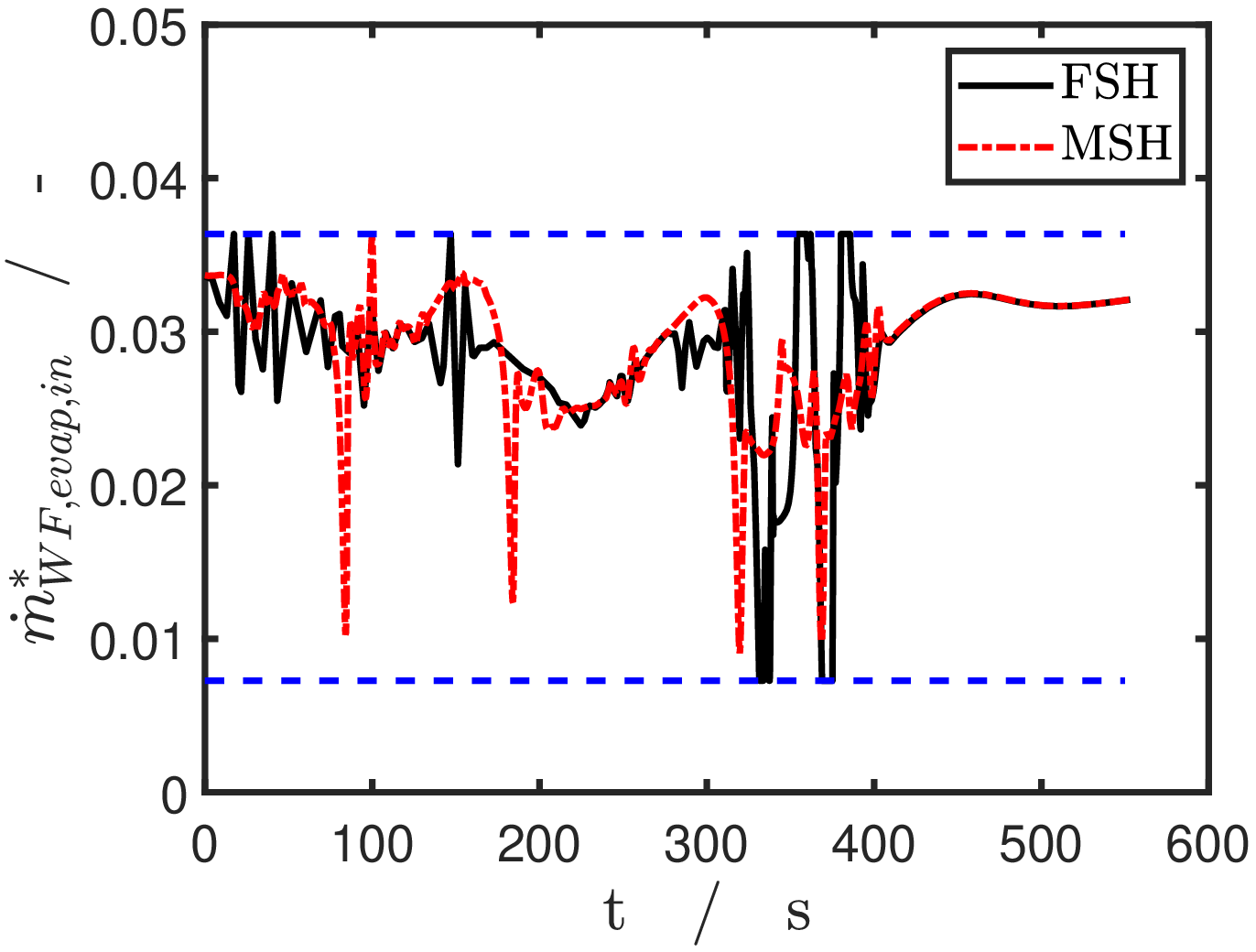}
	\caption{WF mass flow}
	\label{fig:mWfWhtc}
\end{subfigure}
~
\begin{subfigure}[h]{0.45\linewidth}
	\includegraphics[width=\textwidth]{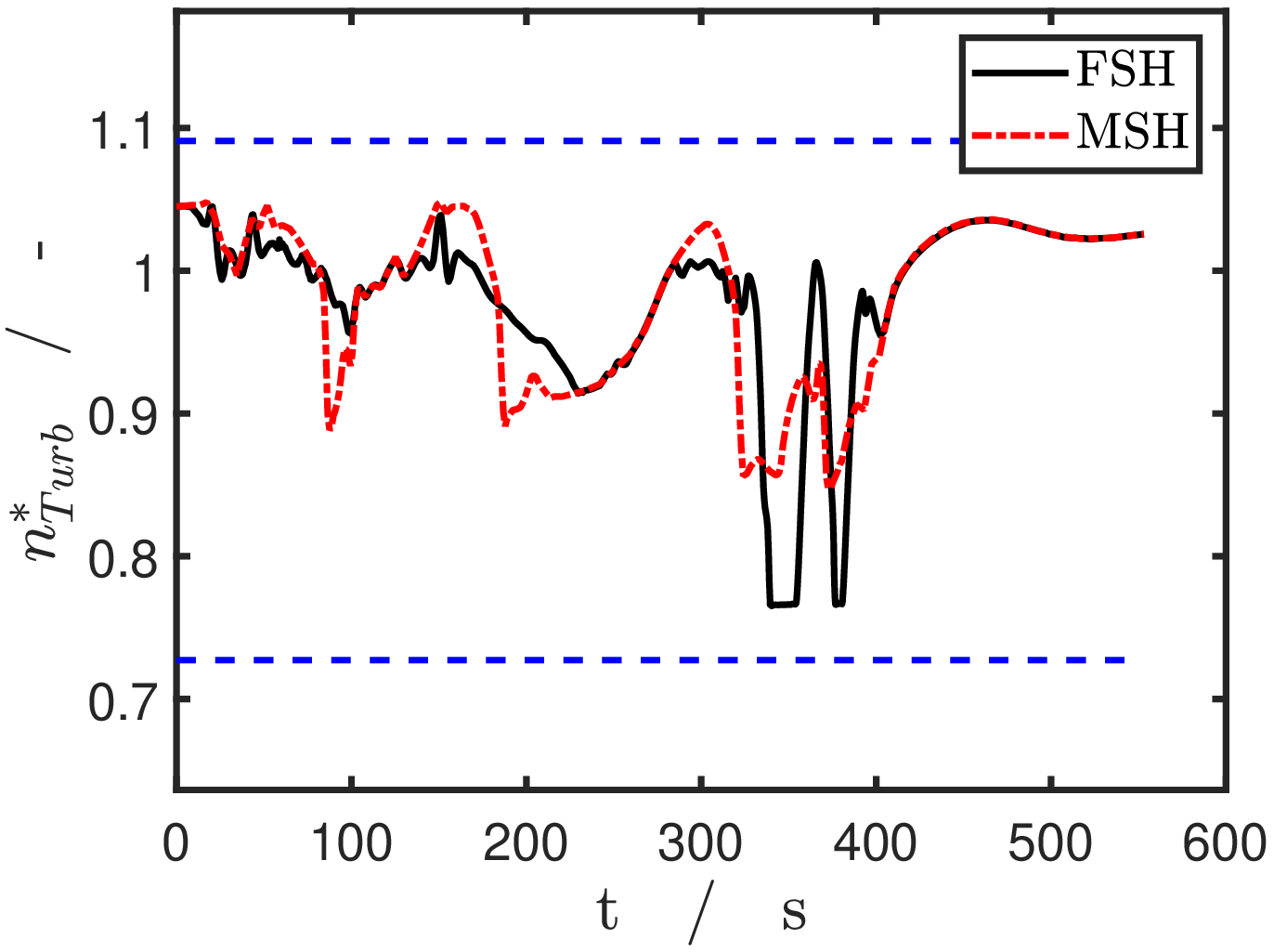}
	\caption{Turbine speed}
	\label{fig:nTurbWhtc}
\end{subfigure}
\\
\begin{subfigure}[h]{0.45\linewidth}
	\includegraphics[width=\textwidth]{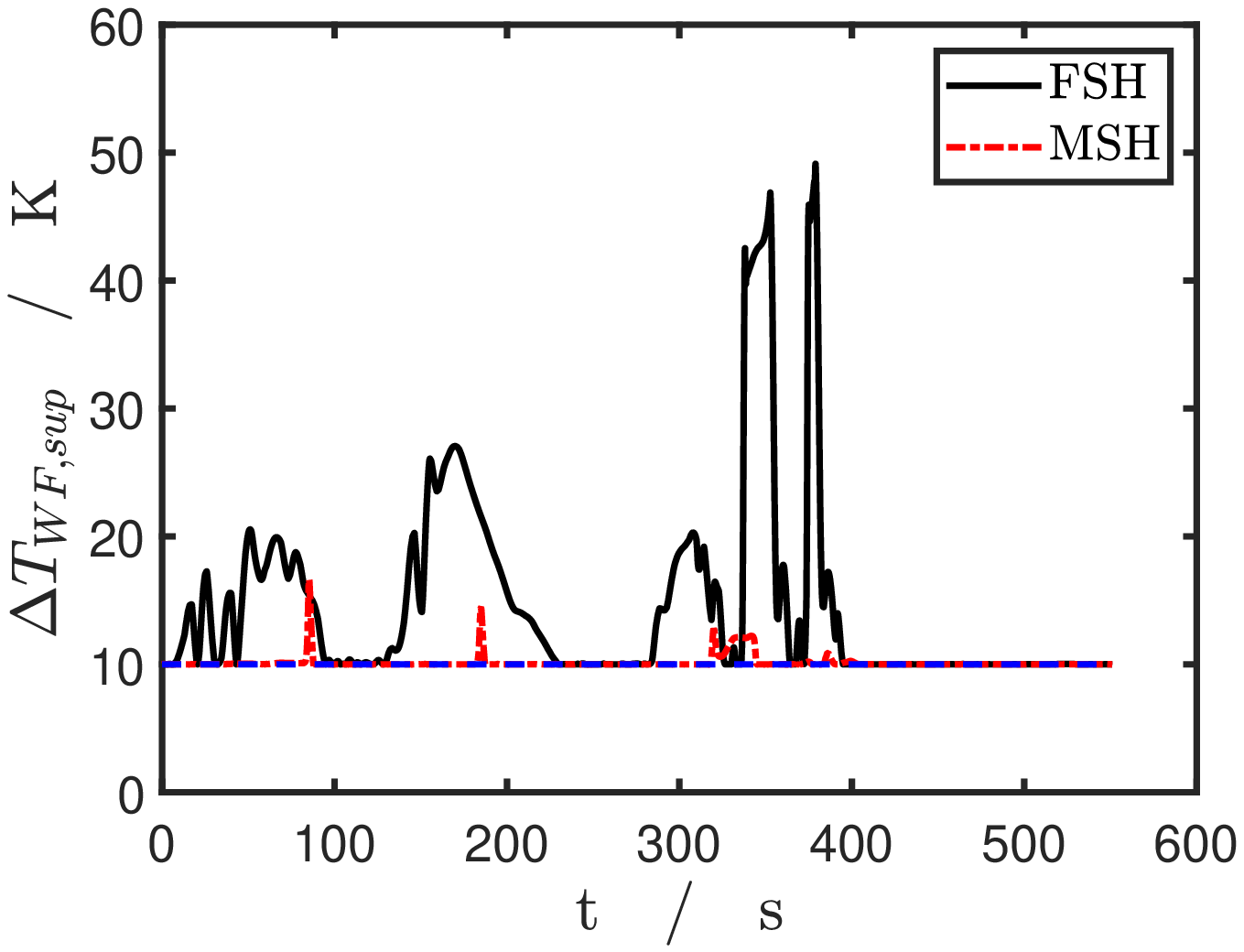}
	\caption{Evaporator superheating}
	\label{fig:tSupWhtc}
\end{subfigure}
~
\begin{subfigure}[h]{0.45\linewidth}
	\includegraphics[width=\textwidth]{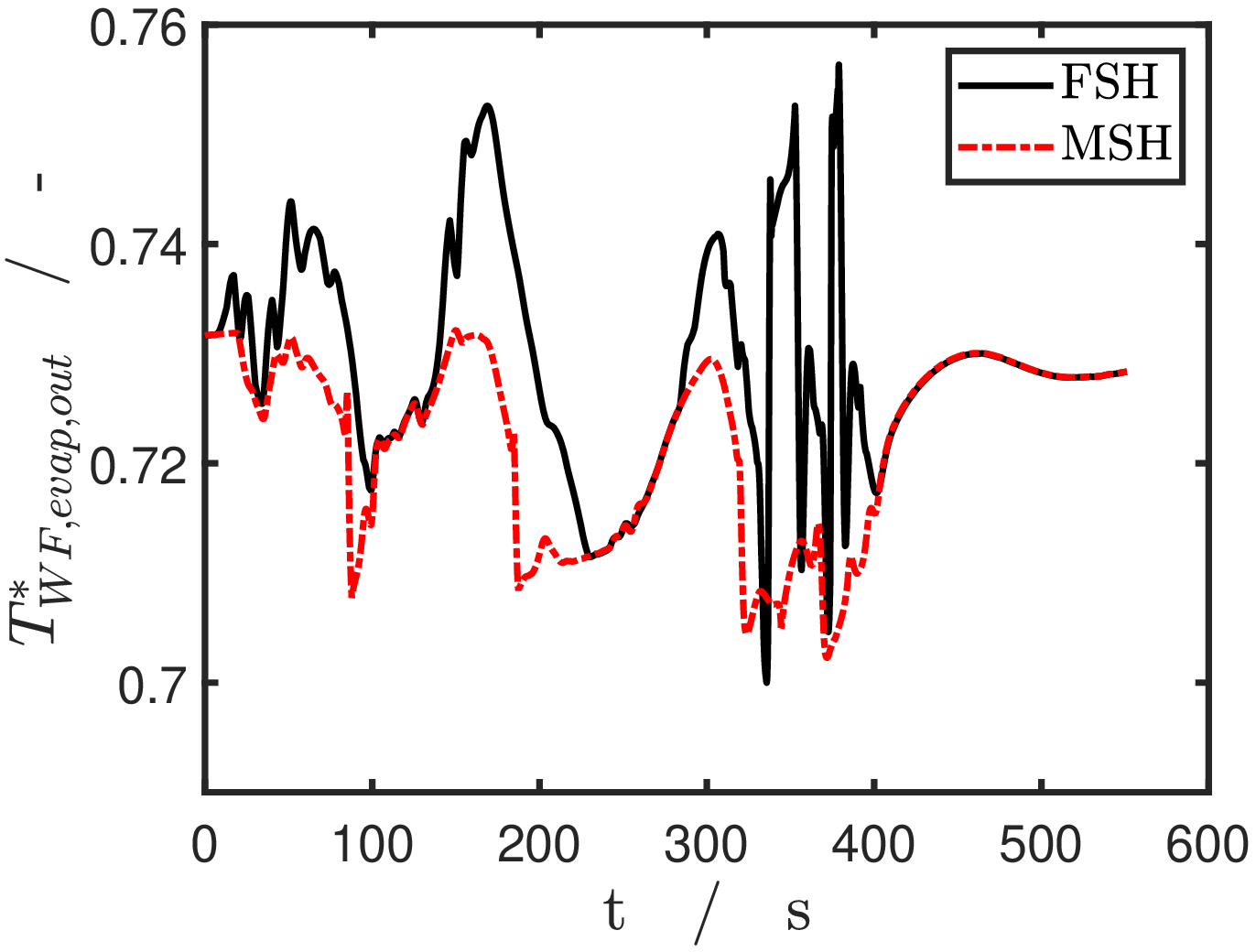}
	\caption{WF evaporator outlet temperature}
	\label{fig:tWfOutWhtc}
\end{subfigure}
\\
\begin{subfigure}[h]{0.45\linewidth}
	\includegraphics[width=\textwidth]{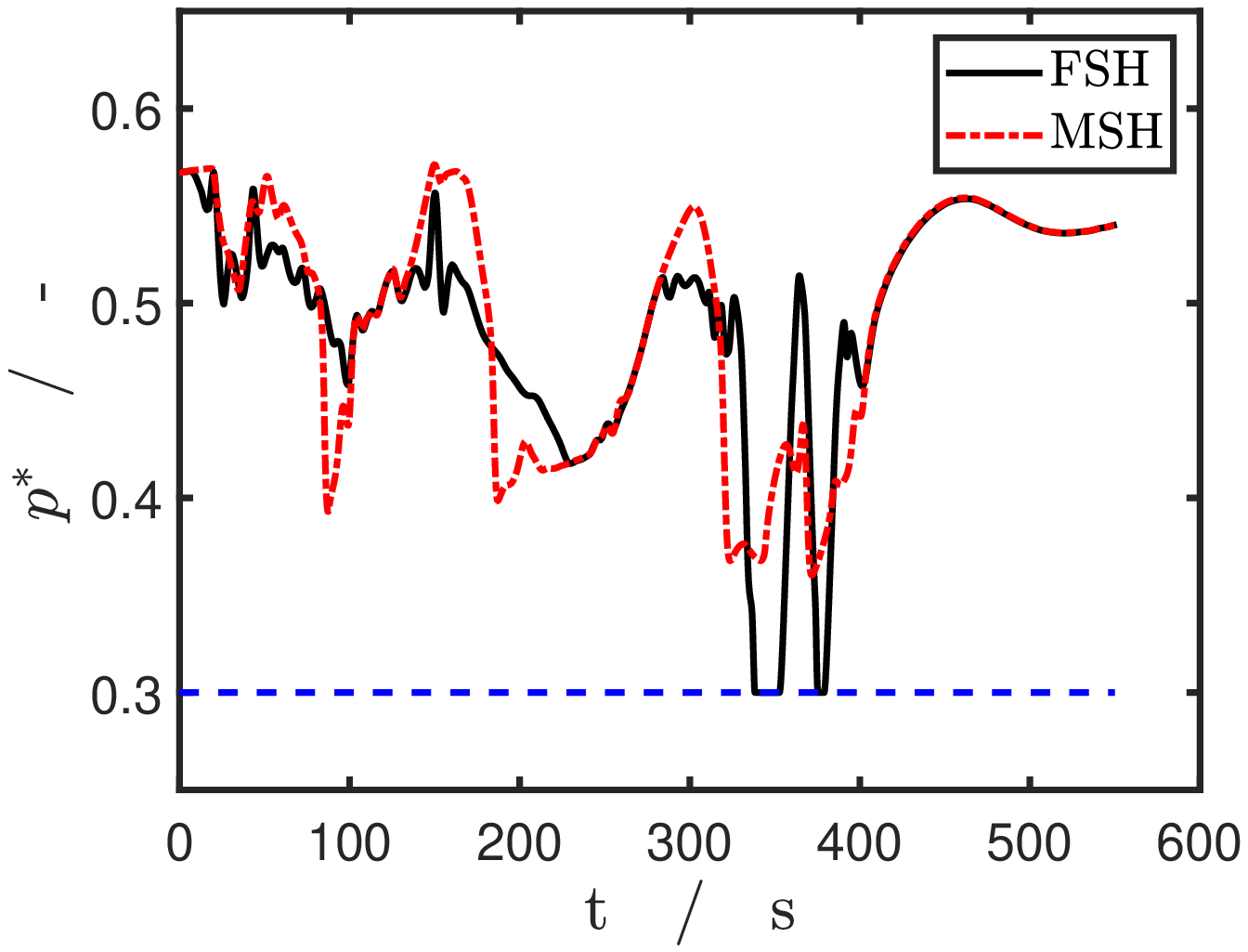}
	\caption{Evaporator pressure}
	\label{fig:pWhtc}
\end{subfigure}
~
\begin{subfigure}[h]{0.45\linewidth}
	\includegraphics[width=\textwidth]{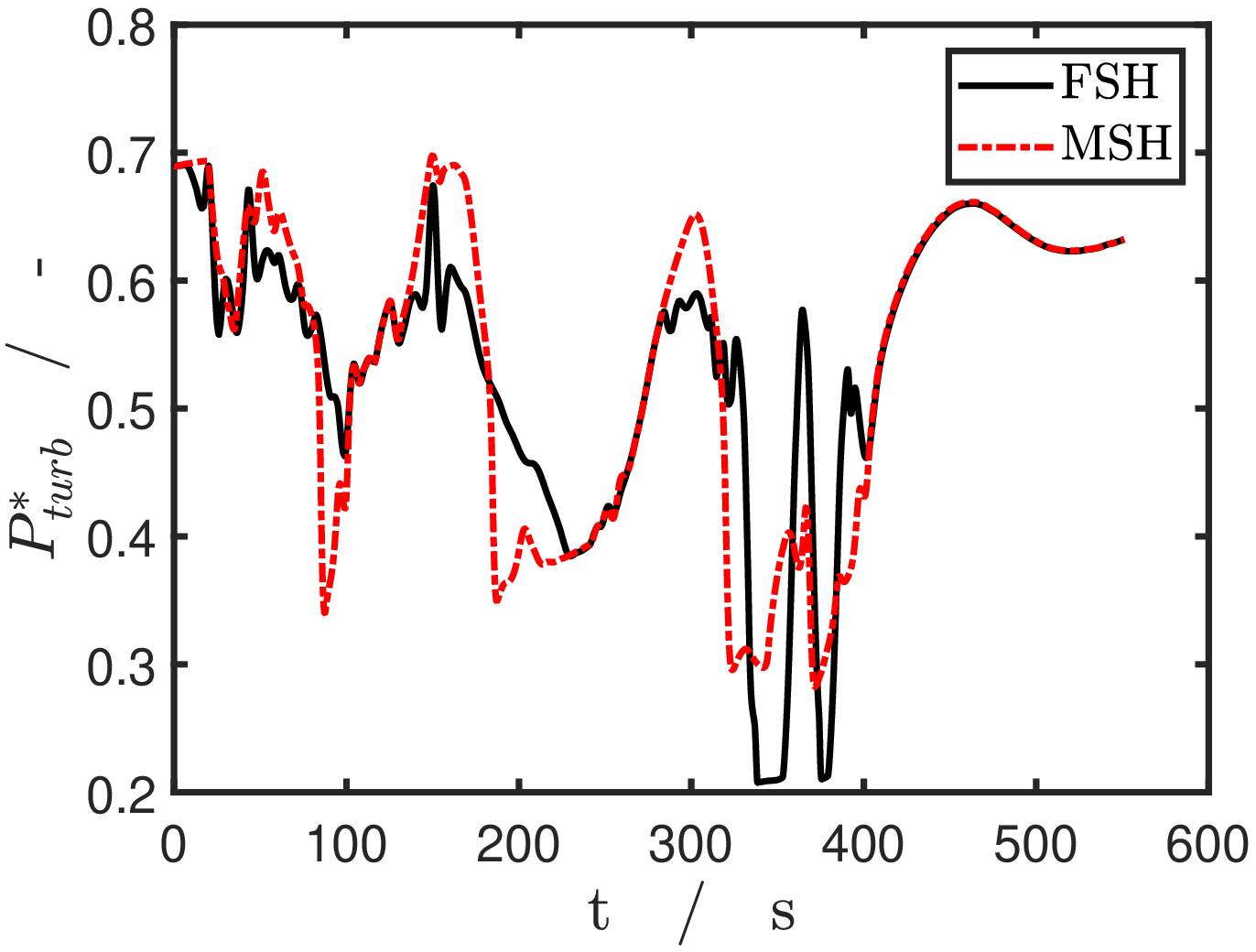}
	\caption{Turbine power}
	\label{fig:pTurbWhtc}
\end{subfigure}
\caption{Optimization results; dashed blue lines indicate lower and upper bounds.}
\label{fig:whtcOpt}
\end{figure}

In this section, we present the optimization of a typical transient driving cycle as it would occur in street traffic.
We assume that there is no limit on expander power at any time and only the safety-related path constraints \eqref{eq:pSafety}-\eqref{eq:supSafety} apply.
The heat source data is taken from an experiment that include parts of the World Harmonized Transient Cycle (Fig.~\ref{fig:heatSourceWhtc}).
For FSH we minimize $\Phi_{2}$.
\\
For MSH, we exploit that $n_{turb}^*$ only optimizes turbine power and use a two-step procedure.
We first minimize $\Phi_{1}$ to obtain minimal superheat and subsequently, we minimize $\Phi_{2}$, where we fix the trajectory of $\dot{m}_{WF,in}^*$ to the optimal solution of the first step and leave $n_{turb}^*$ as DOF to optimize turbine power.
As we can separate the two optimization tasks of achieving minimal superheat and obtaining optimal turbine operation for the resulting operating conditions, we can avoid a weighting between those objectives.
All optimization problems are subject to \eqref{eq:differential}-\eqref{eq:xBpvConstraint}.
\\
\\
The optimized trajectories for the DOF and key variables are depicted in Fig.~\ref{fig:whtcOpt}.
No trajectories for the exhaust bypass valves are presented as it remains fully opened at all times for both policies.
For the examined case, a value of $P_{net,av}^*=0.5307$ is obtained with FSH and $P_{net,av}^*=0.5280$ with MSH.
In other words, the additional  flexibility merely yields a 0.5\% increase in net average power.
However, the trajectories differ strongly from each other in this case and several deviations from minimal superheat occur for FSH while for MSH only small deviations occur in order to maintain feasibility.
Between $t=0~\mathrm{s}$ and $t=100~\mathrm{s}$, the peaks for FSH occur at a high frequency and do not exceed 25~K.
This might be due to the adaption algorithm as it is well known that a very fine discretization can lead to oscillatory control profiles \cite{Schlegel2005}.
More interesting are the two largest deviations in Fig.~\ref{fig:tSupWhtc}, which do occur between $t=300~\mathrm{s}$ and $t=400~\mathrm{s}$ with two peaks exceeding 40~K,  where the exhaust gas exhibits comparatively mild fluctuations.
Apparently, the optimizer exploits the fact that temporarily operating at higher superheat, hence lower pressure level, can be advantageous.
Further analysis showed that the amount of the heat recovered from the exhaust gas and transferred to the WF is higher for MSH.
Moreover, we ruled out that the behavior is due to the fluid-dependent turbine efficiency map by executing the optimization with a turbine with fixed efficiencies where the behavior persisted.
As a test using a different WF (cyclopentane) did not exhibit any peaks in superheat, the behavior appears to be fluid-specific.
Although it is of academic interest, further investigation of the behavior is beyond the scope of this manuscript.
\\
\\
As can be seen from Fig.~\ref{fig:mWfWhtc} and \ref{fig:nTurbWhtc}, the WF mass flow and turbine speed exhibit fast changes for FSH.
It is questionable if a physical unit would be able to follow these trajectories and to what extent the additional strain would results in reduced lifetime of the components.
For MSH, actor action is less drastic and less mechanical strain is expected.
Further, realizing the peaks in superheat in a control setting would require foresight of the exhaust gas conditions.
Considering these observations and the fact that MSH only produces 0.5\% less power than FSH, using minimal superheat seems to be an appropriate control objective, when no other constraints apply.

\section{Optimal operation including limitations on turbine power} \label{sec:activeConstraints}
The case presented in Sec.~\ref{sec:nominal} can be considered as a best case scenario as no constraints beyond the safety constraints are considered.
The strongest assumption we made in Sec.~\ref{sec:nominal} is that the power produced by the turbine can always be utilized completely.
Further, operational constraints, i.e., \eqref{eq:pSafety} or \eqref{eq:tSafety}, could become active.
We consider these scenarios in this section.
\subsection{Active power constraint} \label{sec:activeConstraints1}
\begin{figure}[H]
\centering
\begin{subfigure}[h]{0.45\linewidth}
	\includegraphics[width=\textwidth]{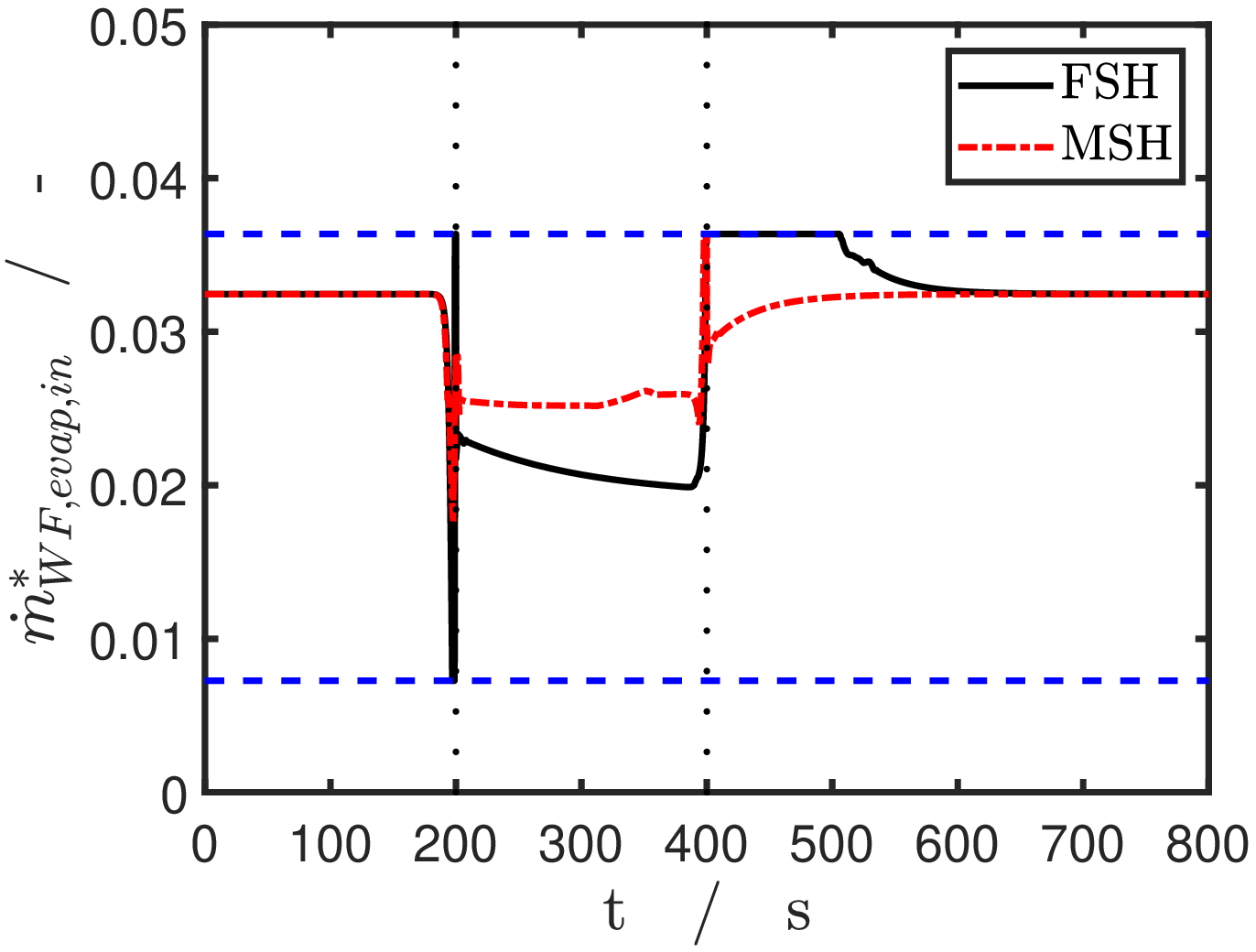}
	\caption{WF mass flow}
	\label{fig:mWfPowLim}
\end{subfigure}
~
\begin{subfigure}[h]{0.45\linewidth}
	\includegraphics[width=\textwidth]{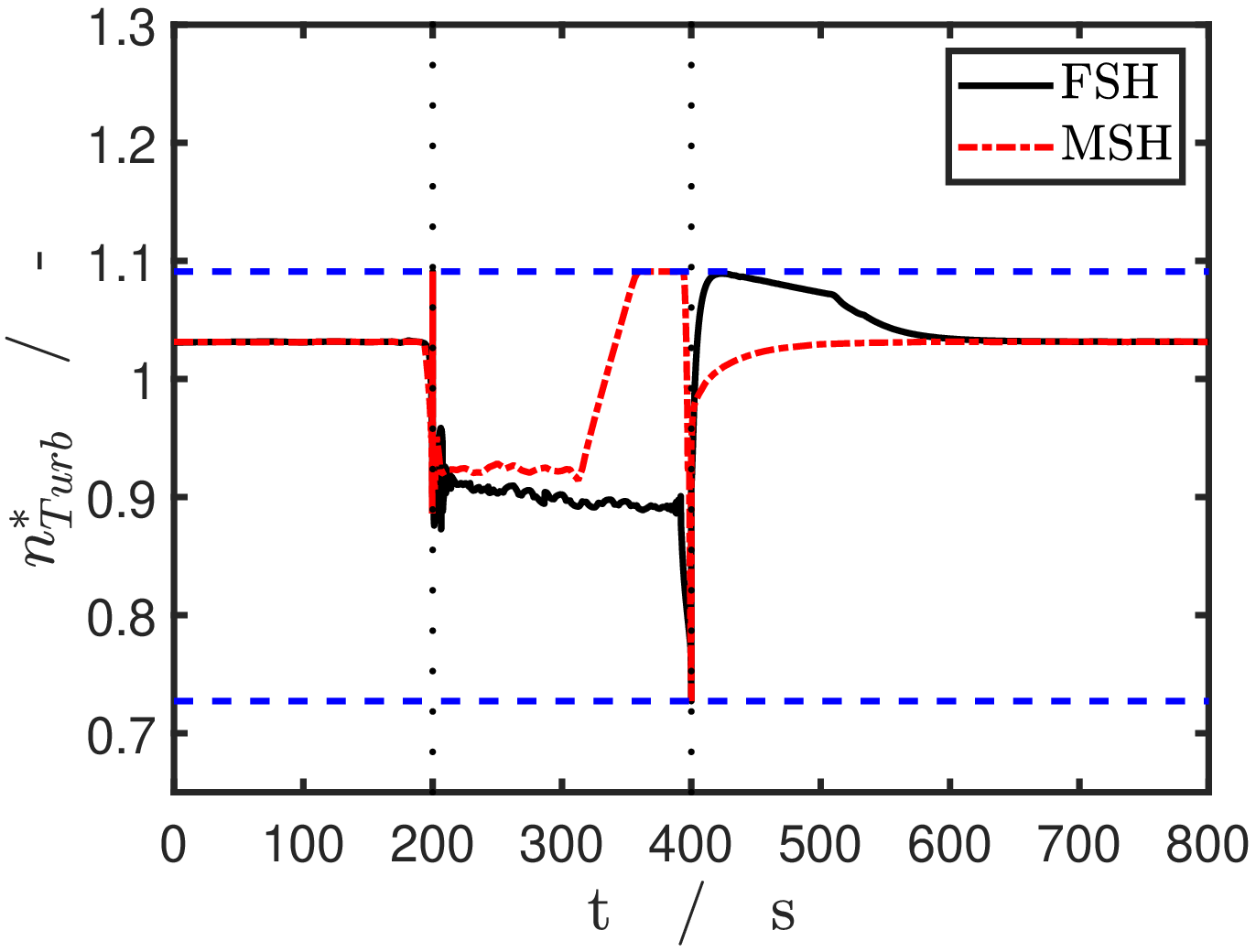}
	\caption{Turbine speed}
	\label{fig:nTurbPowLim}
\end{subfigure}
\\
\begin{subfigure}[h]{0.45\linewidth}
	\includegraphics[width=\textwidth]{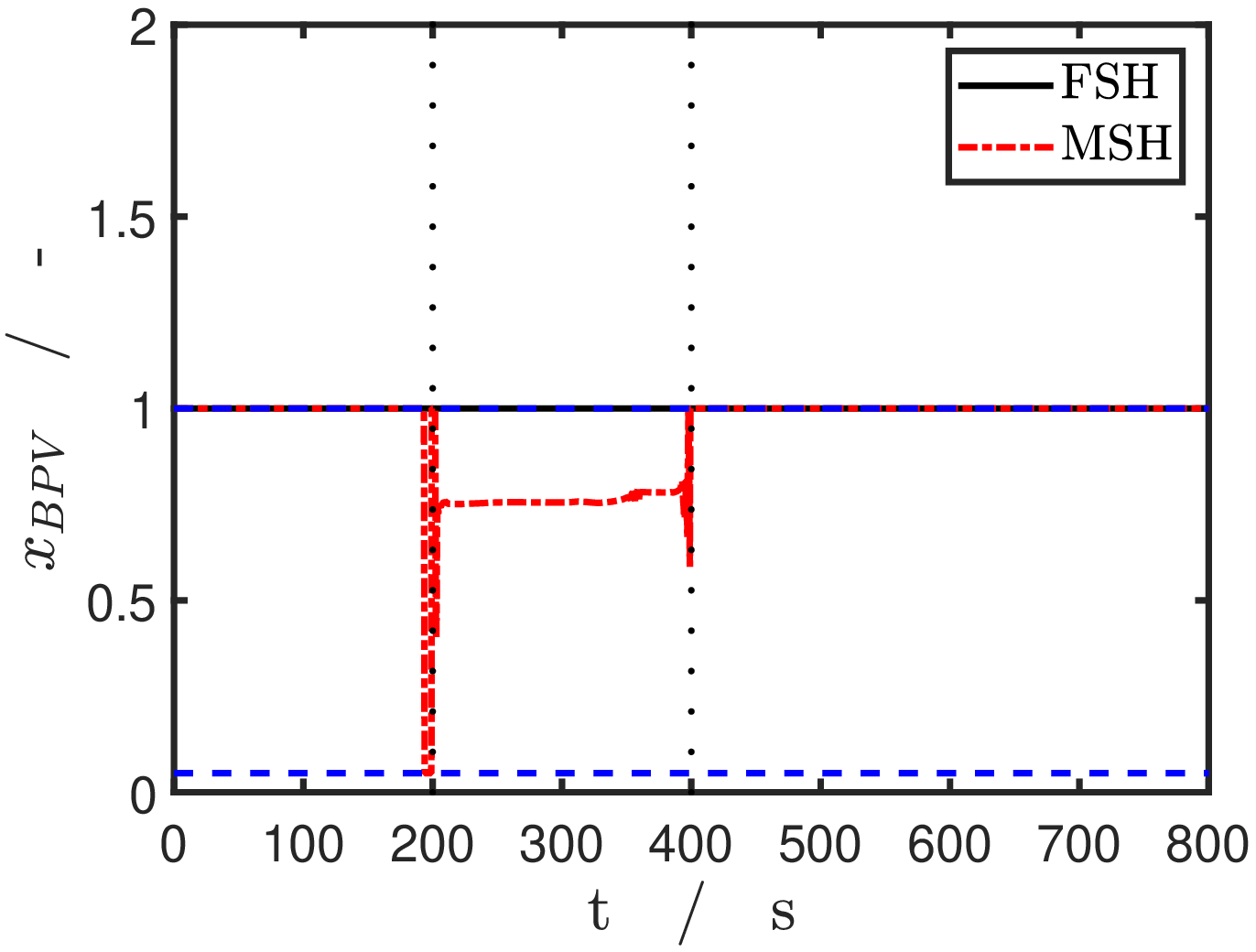}
	\caption{Bypass valve position}
	\label{fig:xBpPowLim}
\end{subfigure}
~
\begin{subfigure}[h]{0.45\linewidth}
	\includegraphics[width=\textwidth]{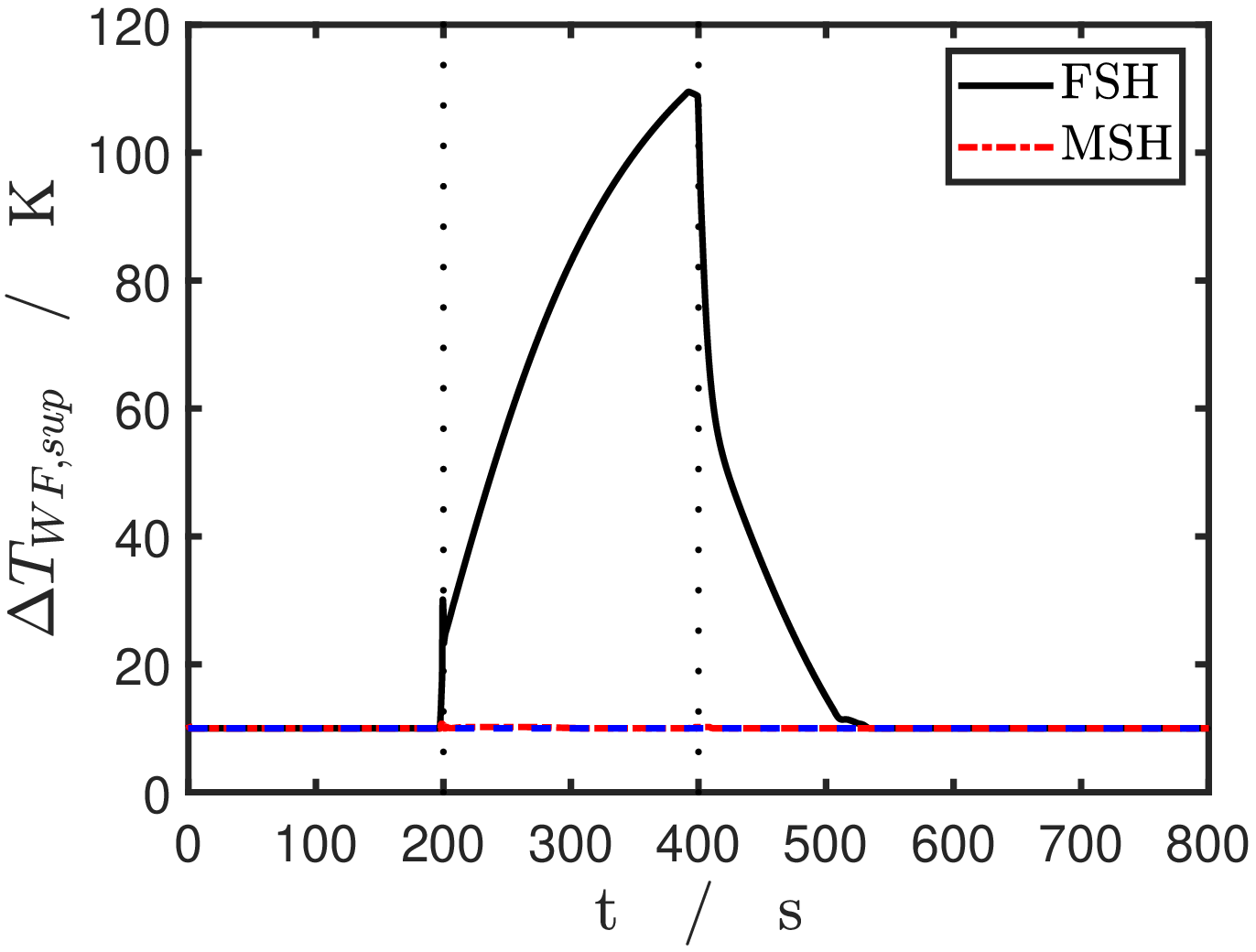}
	\caption{Evaporator superheating}
	\label{fig:tSupPowLim}
\end{subfigure}
\\
\begin{subfigure}[h]{0.45\linewidth}
	\includegraphics[width=\textwidth]{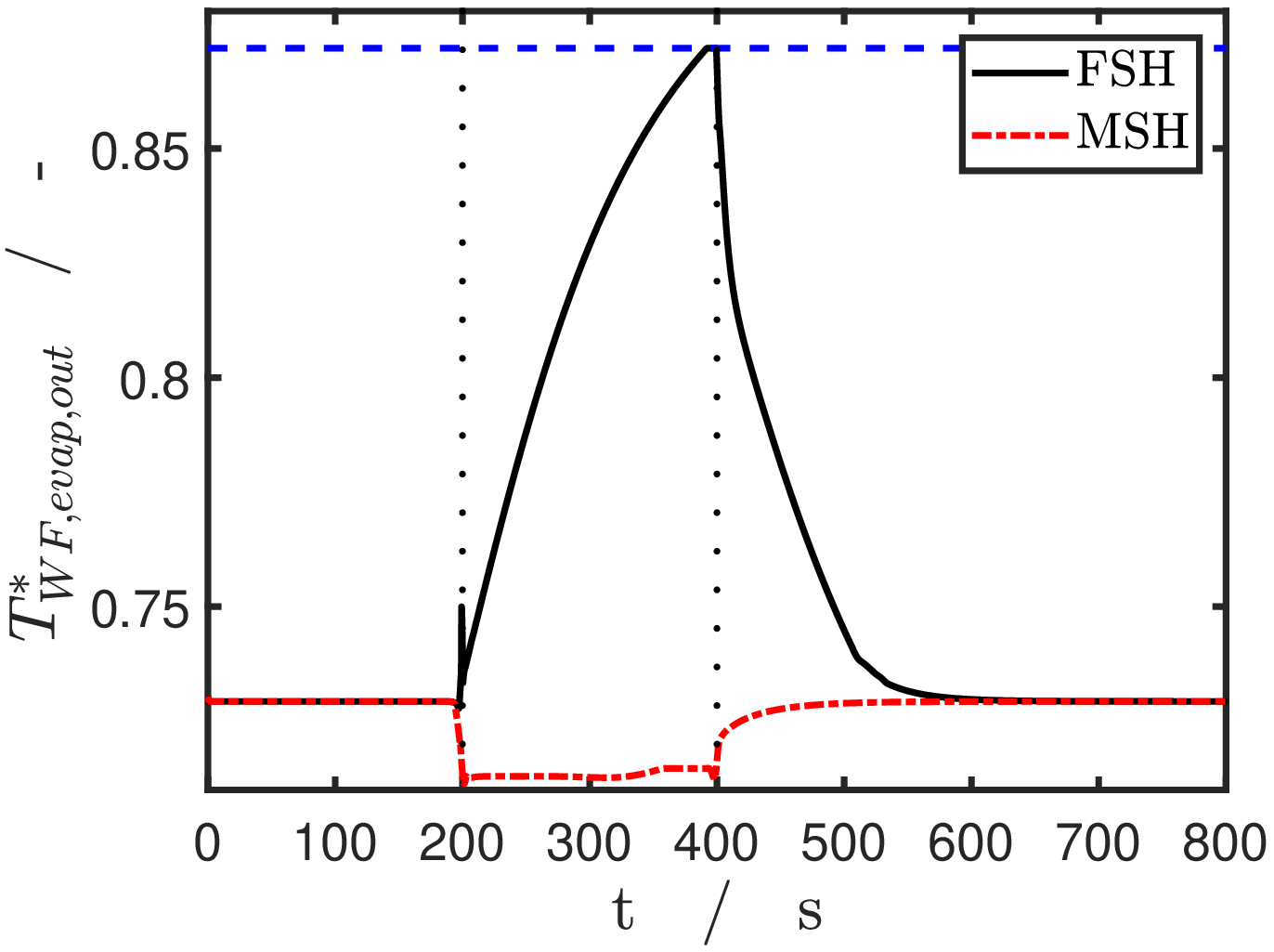}
	\caption{WF evaporator outlet temperature}
	\label{fig:tWfOutPowLim}
\end{subfigure}
~
\begin{subfigure}[h]{0.45\linewidth}
	\includegraphics[width=\textwidth]{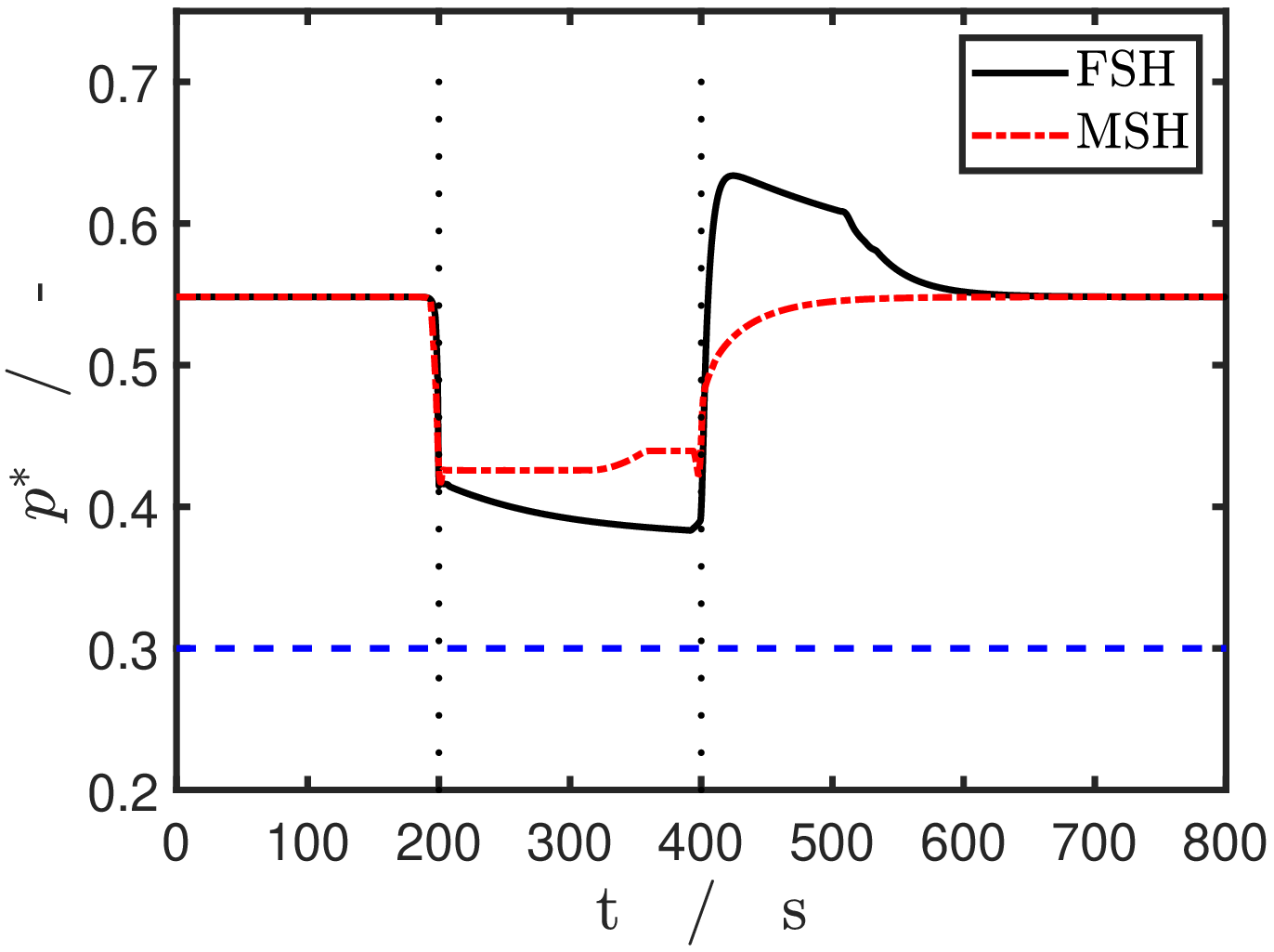}
	\caption{Evaporator pressure}
	\label{fig:pPowLim}
\end{subfigure}
\\
\begin{subfigure}[h]{0.45\linewidth}
	\includegraphics[width=\textwidth]{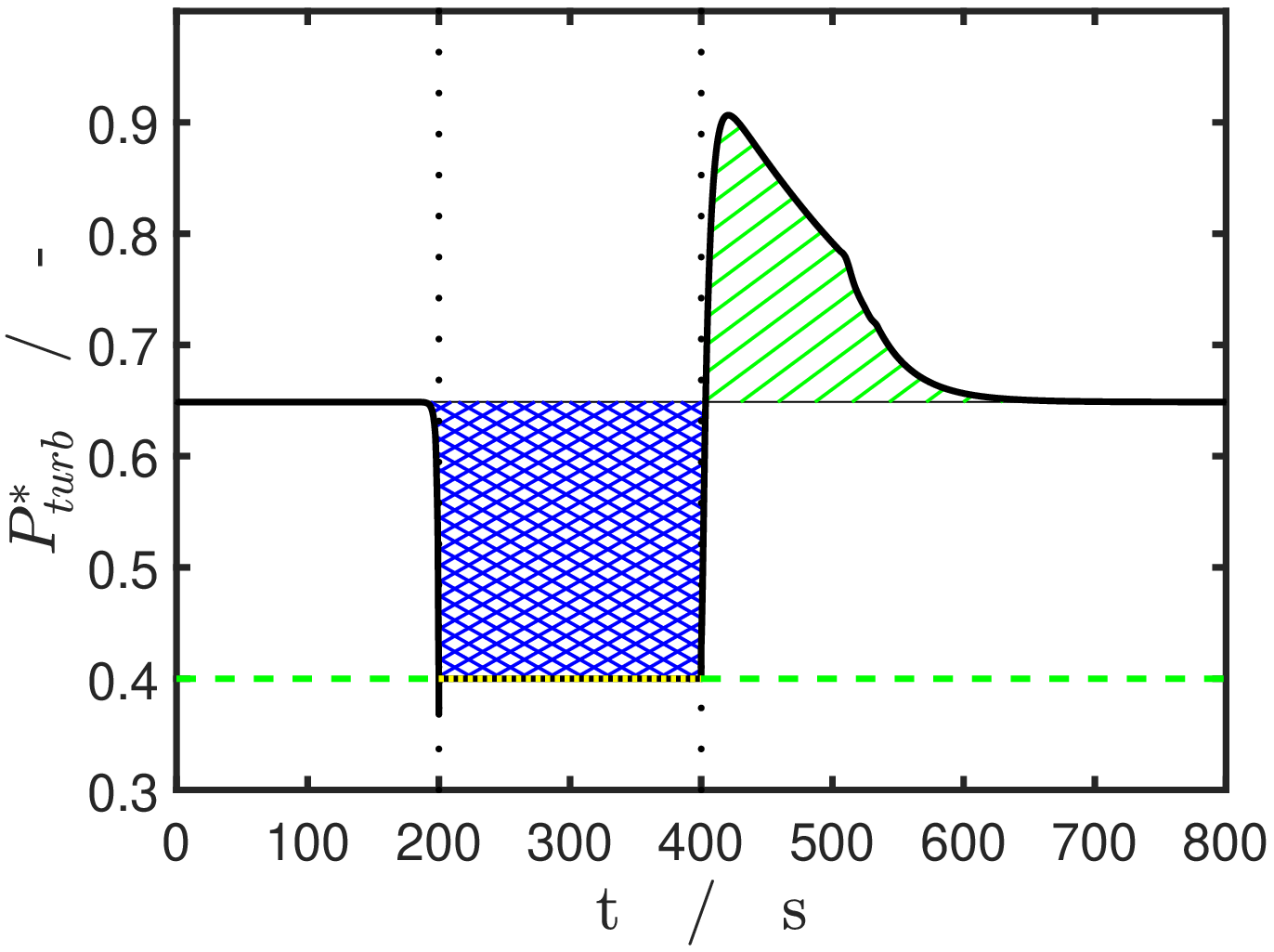}
	\caption{Turbine power for FSH}
	\label{fig:Turbine power}
\end{subfigure}
~
\begin{subfigure}[h]{0.45\linewidth}
	\includegraphics[width=\textwidth]{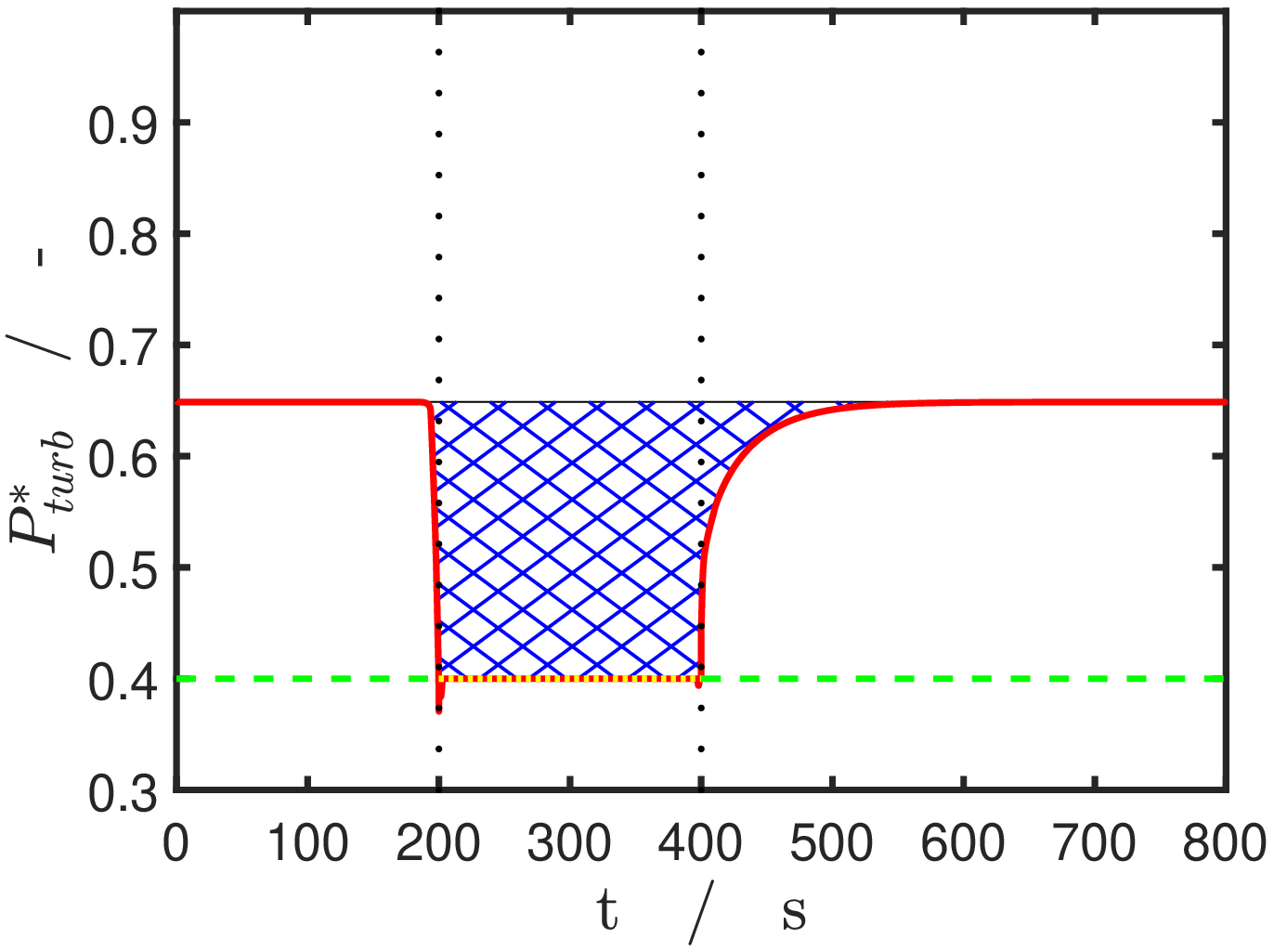}
	\caption{Turbine power for MSH}
	\label{fig:Turbine powerMinSup}
\end{subfigure}
\caption{Results of the optimization for the power limitation case, dashed blue lines indicate lower and upper bounds. Bounds that do not apply at all times are depicted as dashed green lines when they do not apply and yellow dotted lines when they apply.}
\label{fig:trajectoriesPowLim}
\end{figure}

\begin{figure}[h!]
\centering
\includegraphics[width=0.9\linewidth]{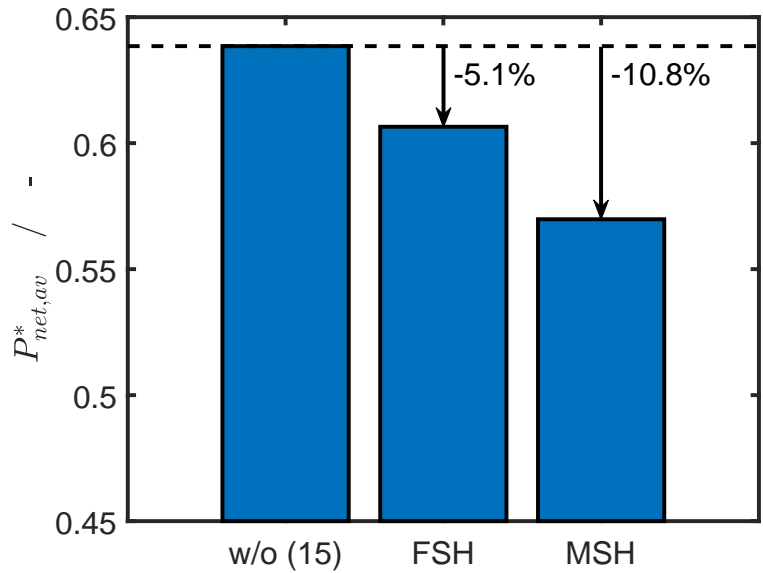}
\caption{Normalized net average power produced for both policies compared to the case where no power limitation occurs (i.e., without constraint \eqref{eq:powLim}). Clearly, FSH exhibits superior performance over MSH.}
\label{fig:tablePowLim}
\end{figure}

To account for a situation where only a limited turbine power can be utilized, e.g., due to maximal charging current of a battery system, we consider a scenario with $t_f=800~\text{s}$ and constant exhaust gas conditions with $\dot{m}_{exh}^*=0.201$, $T_{exh,in}^*= 0.967$.
We assume that the expander power is temporarily limited. 
Here, we arbitrarily choose $t_1=200~\mathrm{s}$ and $t_2=400~\mathrm{s}$ as the start and end of the interval in which the power limitation applies.
We assume knowledge of the times where the power limitation applies as we do with the heat source signals in order to obtain an upper bound on system performance.
We realize this scenario by formulating a multistage dynamic optimization problem consisting of three stages.
For FSH, we use $\Phi_2$ as objective function.
The optimization problem is subject to \eqref{eq:differential}-\eqref{eq:xBpvConstraint} and the turbine power constraint, active in the second stage, is added in~\eqref{eq:powLim}
\begin{align}
& P_{turb}^*\left(t\right)\leq P_{turb}^{*,max} & \forall t\in\left[t_1,\, t_2\right]\, . \label{eq:powLim}
\end{align}
Due to \eqref{eq:powLim}, minimizing superheat and maximizing turbine power are not independent for this scenario and we do not use the two-step strategy from Sec.~\ref{sec:nominal} for MSH.
Instead, we realize MSH by minimizing $\Phi_{2}$ subject to \eqref{eq:differential}-\eqref{eq:xBpvConstraint}, \eqref{eq:powLim} and adding an upper bound on superheat \eqref{eq:powLimTSupMax}:
\begin{align}
T_{sup}^{max}=10.7~\mathrm{K} \, . \label{eq:powLimTSupMax}
\end{align}
Thereby, we find minimal superheat strategy with maximum power production by providing a small range of the permissible superheat with \eqref{eq:supSafety} and \eqref{eq:powLimTSupMax} and minimizing $\Phi_{2}$.
For reference, we consider the case without \eqref{eq:powLim}, i.e., operation at steady-state to allow for an estimate of the energy that is lost due to the power limitation.
\\
\\
The optimal results of the DOF and relevant variables are presented in Fig.~\ref{fig:trajectoriesPowLim}.
Here, the advantages of FSH are evident.
The optimizer exploits the thermal capacity of the evaporator to store thermal energy during the power limitation which is released after the power limitation ends. 
Shortly before the power limitation phase begins, the WF mass flow is reduced (Fig.~\ref{fig:mWfPowLim}), while the exhaust bypass valve remains fully opened (Fig.~\ref{fig:xBpPowLim}).
Consequently, superheat rises and increases up to more than $100$~K at the end of the power limitation phase.
Approximately at $t=390~\mathrm{s}$, $T_{WF,evap,out}^*$ reaches its upper bound and the WF mass flow is increased, which decreases superheat.
The turbine speed is adjusted to a suboptimal point to satisfy the limitation on power output.
When the power limitation ends at $t=400~\mathrm{s}$, the WF mass flow is set to its maximum value which results in a strong increase in pressure and also in power production.
The periods in time where the turbine power is lower than at optimal steady-state without power limitation are indicated by the blue cross-hatched area and the periods in time where turbine power is higher are indicated by the green hatched area in Fig.~\ref{fig:Turbine power}.
The optimizer exploits that  the heat exchanger wall temperature has increased during the power limitation.
Consequently, a higher WF mass flow can be evaporated as can be seen from Fig.~\ref{fig:mWfPowLim}.
Hence, parts of the energy not used earlier can be recovered.
\\
\\
MSH, however, does not exploit this option to save energy, as can be seen in Fig.~\ref{fig:tSupPowLim}.
To allow for the required reduction of turbine power, the exhaust bypass valve opens shortly before the power limitations begins which can be seen in Fig.~\ref{fig:xBpPowLim}.
During the power limitation, it remains partially opened and a part of exhaust gas is bypassed to allow for satisfaction of the superheat path constraint.
The system approaches a steady-state at minimal superheat with a partially opened exhaust bypass valve, hence not making full use of the exhaust gas potential.
At the end of the power limitation, the valve is closed again and all the exhaust gas passes through the evaporator.
In contrast to FSH, however, there is no heat available that can be recovered from the evaporator walls.
Consequently, the system takes some time to reach the initial steady-state which results in some additional loss in a period where FSH exceeds the steady-state turbine power.
\\
\\
The results that are obtained from visual inspection are supported by Fig.~\ref{fig:tablePowLim} where the resulting $P_{net,av}^*$ for each operating policy is presented.
FSH avoids 53\% of the losses associated with MSH.
This result emphasizes that operating at minimal superheat is not necessarily always the best policy.
The behavior can be implemented in a control strategy, albeit imperfectly, without knowledge about the future exhaust conditions.
\subsection{Power limitation with high exhaust gas mass flow}
\begin{figure}[h!]
\centering
\begin{subfigure}[h]{0.45\linewidth}
	\includegraphics[width=\textwidth]{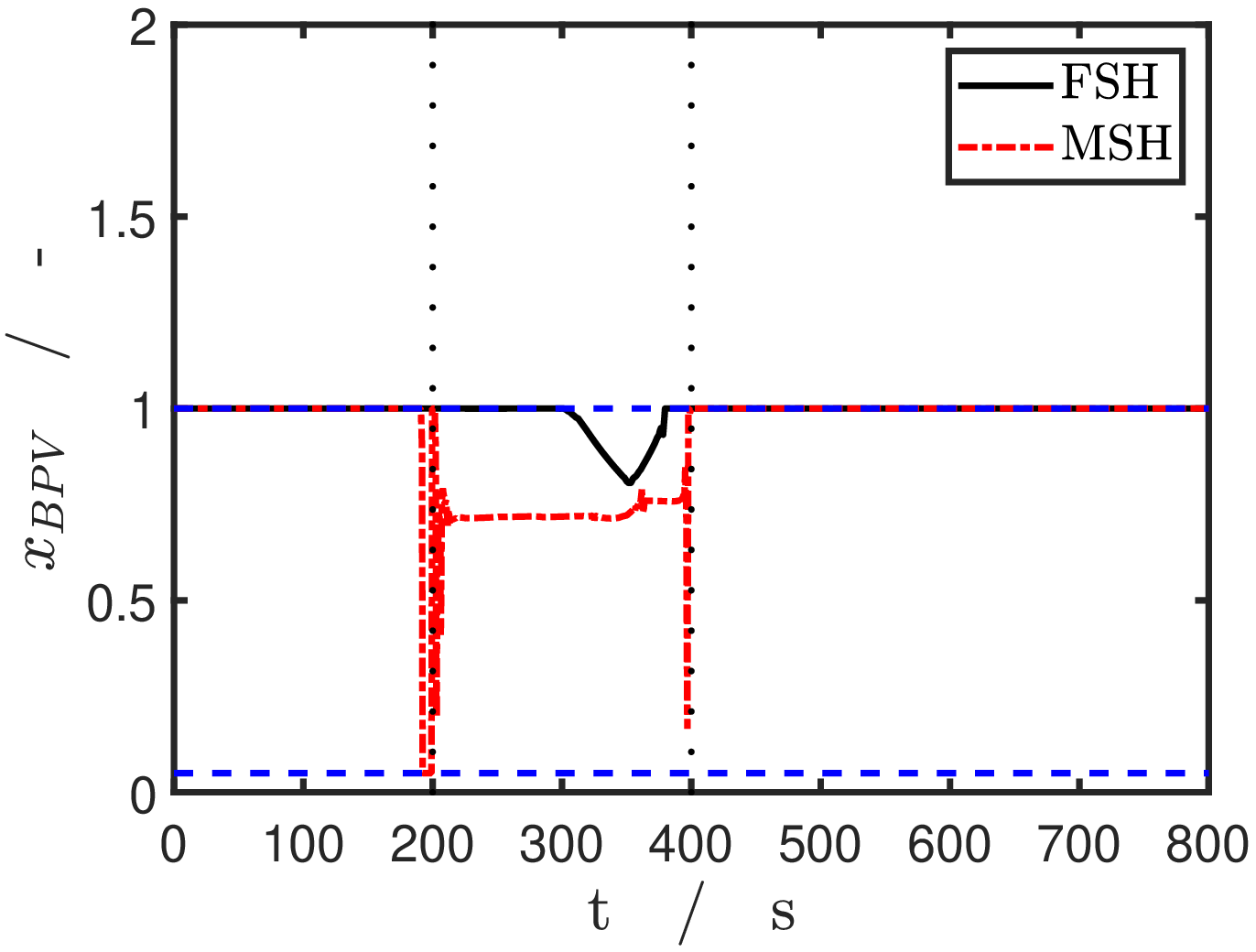}
	\caption{Bypass valve position}
	\label{fig:xBpPowLimExt}
\end{subfigure}
~
\begin{subfigure}[h]{0.45\linewidth}
	\includegraphics[width=\textwidth]{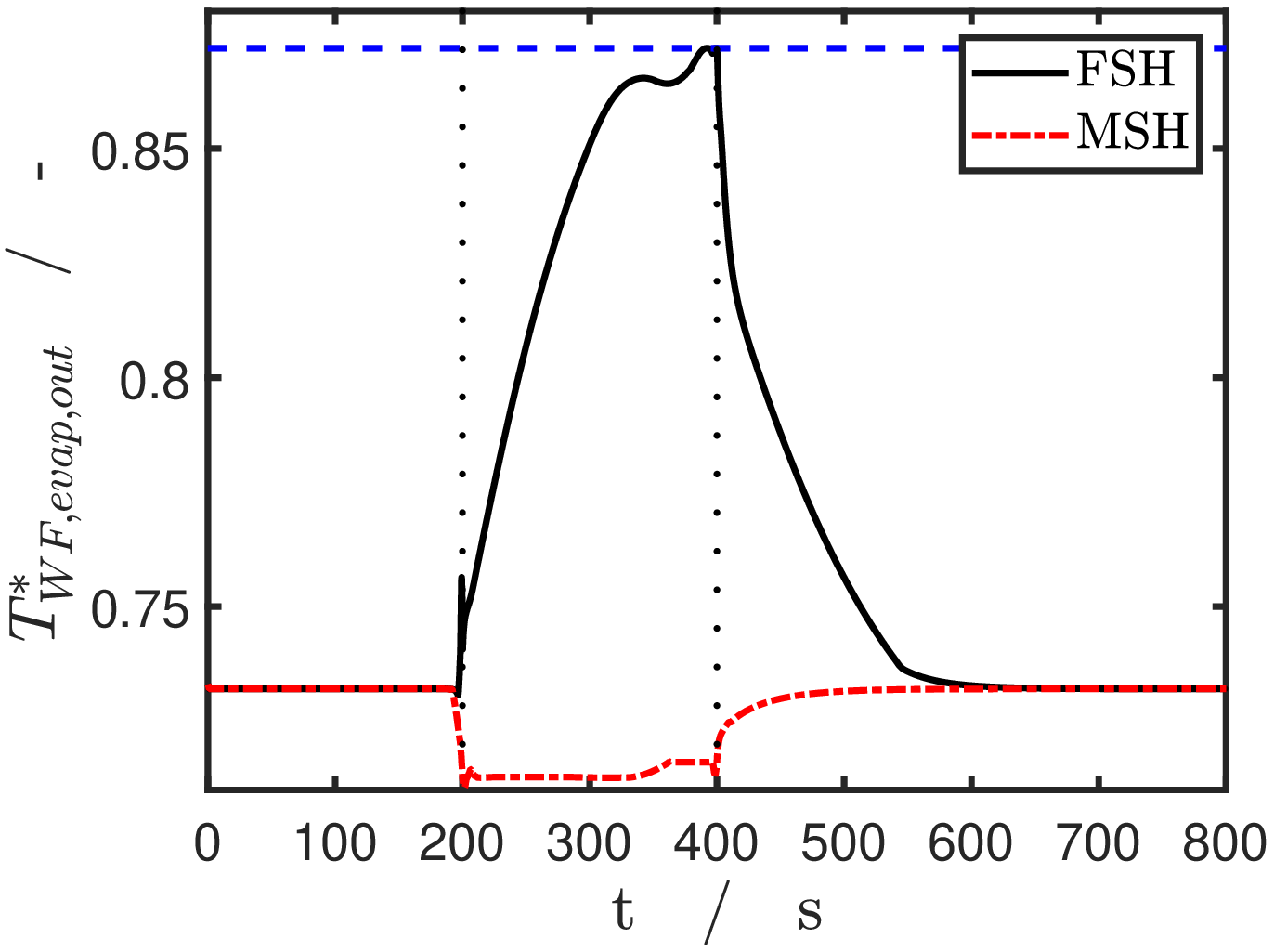}
	\caption{WF evaporator outlet temperature}
	\label{fig:tWfOutPowLimExt}
\end{subfigure}
\caption{Results of the optimization for the power limitation case, dashed blue lines indicate lower and upper bounds.}
\label{fig:trajectoriesPowLimExt}
\end{figure}

\begin{figure}[h!]
\centering
\includegraphics[width=0.9\linewidth]{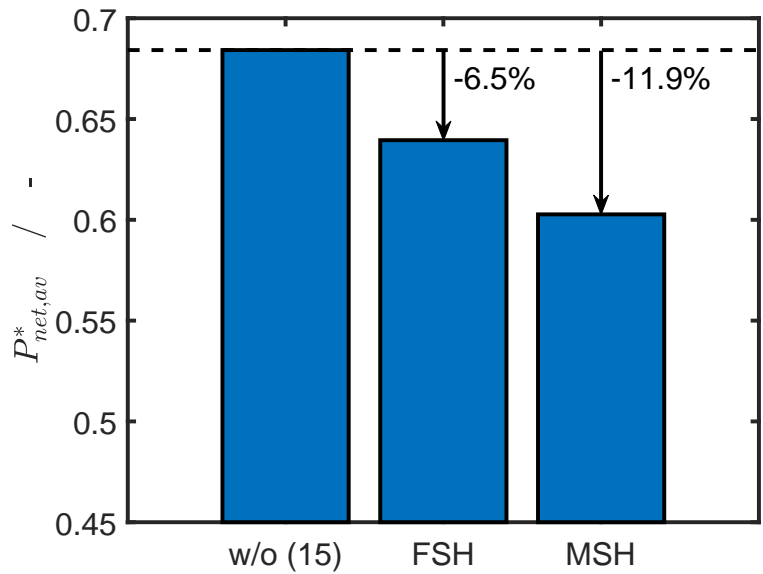}
\caption{Normalized net average power produced for both policies compared to the case where no power limitation occurs (i.e., without constraint \eqref{eq:powLim}) with $\dot{m}_{exh}$ increased by 5\% compared to Fig.~\ref{fig:tablePowLim}. The relative advantage of FSH over MSH shrinks when additional constraints become active.}
\label{fig:tablePowLimExt}
\end{figure}

As can be seen from Fig.~\ref{fig:tWfOutPowLim}, the WF outlet temperature is briefly maintained at its upper bound.
The optimizer is, however, capable of preventing the use of the exhaust bypass valve by increasing the pressure and choosing a suboptimal turbine speed.
It is clear that either a higher exhaust gas mass flow, temperature or longer duration of the power limitation will result in a situation, where it will be required to bypass some of the exhaust gas which will reduce the benefits of FSH.
To assess such a scenario, we increase the exhaust gas mass flow chosen in Sec.~\ref{sec:activeConstraints1} by 5\% to $\dot{m}_{exh}^*=0.211$ and solve the same optimization problems.
\\
\\
The resulting signal of the exhaust bypass valve position and the resulting trajectory of the WF outlet temperature are presented in Fig.~\ref{fig:trajectoriesPowLimExt}.
For FSH, the exhaust gas bypass valve is partially opened between $t\approx 300~\mathrm{s}$ and $t\approx 400~\mathrm{s}$ to avoid the WF from exceeding the temperature limit and a portion of the exhaust gas is bypassed (Fig.~\ref{fig:xBpPowLimExt}).
It should be noted that the exhaust bypass valve is reopened before $t=400~\mathrm{s}$ in anticipation of the end of the power limitation so that the WF temperature reaches its upper bound exactly at that point in time (Fig.~\ref{fig:tWfOutPowLimExt}).
This behavior requires a-priori knowledge of the exhaust gas profile and cannot be directly included in the control strategy.
Rather, the valve would be opened once the power limitation would end.
The losses associated with this fact, however, should be negligible and further assessment is beyond the scope of this manuscript.
\\
\\
For MSH, the qualitative behavior is similar to Sec.~\ref{sec:activeConstraints1}.
As FSH requires bypassing a portion of the exhaust gas for the considered scenario, the avoided losses are smaller than in the previous case study, as can be seen from Fig.~\ref{fig:tablePowLimExt}.
Here, only 45\% of the losses associated with MSH can be avoided.
For higher exhaust gas mass flows, FSH consequently results in a reduced relative advantage.
This would also apply for longer power limitations or higher exhaust gas temperatures.
The results from this section suggest that the exhaust bypass valve will only be required for control as a manipulated variable to maintain safe operation.
A general quantitative statement on the benefits of this strategy cannot be made here as it clearly depends on the system at hand and its operating conditions.

\section{Implications on control strategy} \label{sec:discussion}
In Sec.~\ref{sec:nominal} we found that economically optimal dynamic operation exhibits occasional peaks in superheating for a highly transient exhaust gas profile and ethanol as WF.
The gain in produced power, however, is negligible in comparison to operation at minimal superheat.
As the knowledge of the exhaust gas conditions is unrealistic but required to exploit the effect and the resulting control action put unnecessary strain on the actors, it is appropriate to operate the system at minimal superheat during regular operation.
To obtain an optimal turbine speed, a separate optimization problem has to be solved.
However, this can be approximated with a cheap steady state optimization.
The implications of further restrictions beyond safety constraints on the control strategy are more severe.
In Sec.~\ref{sec:activeConstraints} we found that using a flexible superheat operating policy can be significantly more efficient than a minimal superheat operating policy when the turbine power output is temporarily limited.
This result does not depend on a priori knowledge and can be implemented in control strategies.
It implies to track a turbine power set-point by adjusting WF mass flow and turbine rotational speed.
To achieve this, a steady-state optimization could be carried out that aims at satisfying the power constraint while minimizing WF mass flow to guarantee maximal superheat.
The exhaust bypass valve should only be used to guarantee satisfaction of safety constraints once another constraint becomes active.
The system should then be operated at the steady-state which satisfies this additional constraint.

\section{Conclusion and outlook} \label{sec:conclusion}
We assess the optimal operation of an ORC system for waste heat recovery in a heavy-duty diesel truck.
We obtain optimal trajectories for the DOF by means of dynamic optimization with the open-source software tool DyOS \cite{Caspari2019}.
We compare an operating policy that maximizes the net work (FSH) with a policy that maintains minimal superheat while maximizing turbine power (MSH).
\\
\\
First, we assess optimal operation of the WHR system in a transient driving cycle.
Results obtained with FSH indicate that, most of the time, it is best to operate the system at minimal superheat, which is in agreement with literature for steady-state operation.
However, peaks in superheat do occur but gains in power compared to MSH are negligible.
Further, the occurrence of superheat peaks appears to be fluid-specific.
Hence, we recommend operating the system at minimal superheat during nominal operating mode.
This notion is reflected in many published studies on ORC control.
\\
\\
When further limitations apply, MSH can be suboptimal.
This is illustrated for the case of a limitation in permissible turbine power.
For the scenario where we assume constant exhaust gas mass flow and temperature and a temporary constraint on the turbine power, FSH reveals that during that period, increased superheat is greatly beneficial.
The optimizer exploits that thermal energy which cannot be utilized during the power limitation can be stored in the heat exchanger wall for later use.
In the examined case, this avoids 53\% of the power losses resulting from MSH.
Further investigations reveal that the relative advantage of FSH shrinks in scenarios where other constraints become active.
Eventually, the exhaust bypass valve has to be opened and part of the exhaust gas cannot be used as otherwise the maximum WF temperature would be exceeded.
For the examined case, the avoided losses drop to 45\%.
A similar effect is expected to apply when the duration of the power limitation is increased.
In contrast to our previous work \cite{Ghasemi2013}, the behavior observed in this work is due to dynamic effects.
The optimizer exploits that by temporarily storing thermal energy through increased superheat, more power can be produced overall.
\\
\\
Future work should consider other typical situations, e.g., constraints on the cooling capacity.
This would require to include a condenser model.
Further, start-up situations or situations where the WF cannot be fully evaporated are of interest.
Including them would require a discrete-continuous hybrid heat exchanger model \cite{Vaupel2019} which could be modeled for optimization as suggested in \cite{Caspari2020}.
Finally, transferring the insight gained through dynamic optimization into a feasible control strategy is an important task.

\section*{Acknowledgments}
The work leading to this contribution was funded by the Federal Ministry for Economic Affairs and Energy (BMWi) according to a resolution passed by the German Federal Parliament under grant number 19U14010C.
The authors gratefully acknowledge the financial support of the Kopernikus project SynErgie by the Bundesministerium f\"ur Bildung und Forschung (BMBF) and the project supervision by the project management organization Projekttr\"ager J\"ulich (PtJ).

\appendix
\section{Selected model equations}
Here, we present a description of selected model equations, taken from \cite{Huster2018}, required for understanding the model.
For a full description, including parameter values resulting from a dynamic parameter estimation, we refer the reader to \cite{Huster2018}.
\subsection{Evaporator moving boundary model}
For control volumes with single-phase flow, we get the following mass \eqref{chap3_eq:final_mass_balance} and energy \eqref{chap3_eq:final_energy_balance} balances
\begin{align}
A\,\left(\left(z_a-z_b\right)\,\frac{\text{d}\bar{\rho}}{\text{d}t}+\bar{\rho}\,\frac{\text{d}\left(z_b-z_a\right)}{\text{d}t}\right)+\rho_aA\,\frac{\text{d}z_a}{\text{d}t}-\rho_bA\,\frac{\text{d}z_b}{\text{d}t}=\dot{m}_a-\dot{m}_b ,
\label{chap3_eq:final_mass_balance}
\end{align}
\begin{align}
A\,\left(\left(z_b-z_a\right)\bar{\rho}\frac{\text{d}\bar{h}}{\text{d}t}+\left(z_b-z_a\right)\bar{h}\frac{\text{d}\bar{\rho}}{\text{d}t}+\bar{\rho}\bar{h}\frac{\text{d}\left(z_b-z_a\right)}{\text{d}t}\right)-A\,\left(z_b-z_a\right)\frac{\text{d}p}{\text{d}t}  \nonumber \\
+\rho_ah_aA\frac{\text{d}z_a}{\text{d}t}-\rho_bh_bA\,\frac{\text{d}z_b}{\text{d}t}=\dot{m}_a h_a-\dot{m}_b h_b + b_{WF}\alpha_{WF} \left(z_b-z_a\right) \left(T_{w}-\bar{T}\right),
\label{chap3_eq:final_energy_balance}
\end{align}
where $A$ is the cross-sectional area of the fluid channel and $z$ is the longitudinal coordinate. $\rho$, $T$ and $\dot{m}$ are density, temperature and mass flow of the WF, where the subscripts $a$ and $b$ indicate quantities of the left-hand and right-hand boundary of the zones and the overline indicates averaged quantities.
$t$ is the time, $b_{WF}$ the width of the fluid channel and $\alpha_{WF}$ is the heat transfer coefficient from WF to the wall.
The last term on the right hand side of \eqref{chap3_eq:final_energy_balance} is the heat flow from the wall into the WF.
As $\bar{\rho}$ and $\bar{h}$ are algebraic quantities, we account for their time dependence by constructing their total differential w.r.t. the differential quantities \eqref{chap3_eq:rho_total_deriv} and \eqref{chap3_eq:h_total_deriv} and inserting them into \eqref{chap3_eq:final_mass_balance} and \eqref{chap3_eq:final_energy_balance}, in order to obtain formulations, where only actual differential quantities appear in time derivatives.
\begin{align}
\frac{\text{d}\bar{\rho}}{\text{d}t}&=\frac{\partial \bar{\rho}}{\partial p}\frac{\text{d}p}{\text{d}t}+\frac{\partial \bar{\rho}}{\partial \bar{h}}\frac{\text{d}\bar{h}}{\text{d}t}
\label{chap3_eq:rho_total_deriv}
\end{align}
\begin{align}
\frac{\text{d}\bar{h}}{\text{d}t}&=\frac{1}{2}\left(\frac{\text{d}h_a}{\text{d}t}+\frac{\text{d}h_b}{\text{d}t}\right)
\label{chap3_eq:h_total_deriv}
\end{align}
For the two-phase zone, the mass \eqref{chap3_eq:MBTP} and energy \eqref{chap3_eq:EBTP} balances are
\begin{align}
A\Big(\left(\bar{\gamma}\rho'' + \left(1 - \bar{\gamma}\right)\rho'\right)\frac{\text{d}\left(z_b-z_a\right)}{\text{d}t} + \left(z_b-z_a\right)\Big(\frac{\text{d}\bar{\gamma}}{\text{d}t}\left(\rho'' - \rho'\right)
\nonumber\\
+ \bar{\gamma}\frac{\partial \rho''}{\partial p}\frac{\text{d}p}{\text{d}t} + \left(1 - \bar{\gamma}\right)\frac{\partial \rho'}{\partial p}\frac{\text{d}p}{\text{d}t}\Big)\Big) + \rho_a A \frac{\text{d} z_a}{\text{d}t} - \rho_b A \frac{\text{d} z_b}{\text{d}t} 
= \dot{m}_a - \dot{m}_b ,
\label{chap3_eq:MBTP}
\end{align}
\begin{align}
A\Big(\frac{\text{d} \left(z_b-z_a\right)}{\text{d}t}\left(\bar{\gamma} \rho'' h'' + \left(1 - \bar{\gamma}\right) \rho' h' \right) + \left(z_b-z_a\right) \Big(\frac{\text{d}\bar{\gamma}}{\text{d}t}\left(\rho'' h'' - \rho' h'\right)
\nonumber \\
+ 
\bar{\gamma} h'' \frac{\text{d} \rho''}{\text{d} p} \frac{\text{d} p}{\text{d} t} + \left(1 - \bar{\gamma}\right) h' \frac{\text{d} \rho'}{\text{d} p} \frac{\text{d} p}{\text{d} t} + \bar{\gamma} \rho'' 
\frac{\text{d} h''}{\text{d} p} \frac{\text{d} p}{\text{d} t}+ \left(1 - \bar{\gamma}\right) \rho' \frac{\text{d} h'}{\text{d} p} \frac{\text{d} p}{\text{d} t} \Big)\Big)
\nonumber \\
- A \left(z_b-z_a\right) \frac{\text{d} p}{\text{d} t} + A \rho_a h_a \frac{\text{d}z_a}{\text{d} t}
- A \rho_b h_b \frac{\text{d}z_b}{\text{d}t} \nonumber \\
= \dot{m}_a h_a - \dot{m}_b h_b + b_{WF}\alpha_{WF} \left(z_b-z_a\right) \left(T_{w}-\bar{T}\right) ,
\label{chap3_eq:EBTP}
\end{align}
where $\bar{\gamma}$ is the average void fraction calculated with \eqref{chap3_eq:gamma_bar} and the superscripts $'$ and $''$ indicate quantities at liquid and vapor saturation respectively.
The time derivative of the average void fraction can be expressed by constructing the total differential w.r.t. the differential quantities \eqref{chap3_eq:gamma_partial_t}.
\begin{align}
\bar{\gamma} = \frac{\rho'}{\left(h_0 - h_2\right)\left(\rho' - \rho''\right)^2} \bigg\{\left(h_0 - h_2\right) \rho' + \rho'' \bigg[h_2 - h_0
\nonumber \\
+ \left(h' - h''\right) \ln \left(\frac{ \rho'' \left(h'' - h_0\right)}{\rho' \left(h_2 - h'\right)}\right)\bigg]\bigg\}
\label{chap3_eq:gamma_bar}
\end{align}
\begin{align}
\frac{\text{d} \bar{\gamma}}{\text{d} t} = \frac{\partial \bar{\gamma}}{ \partial h_0} \frac{\text{d} h_0}{\text{d} t} + \frac{\partial \bar{\gamma}}{ \partial h_2} \frac{\text{d} h_2}{\text{d} t} + \frac{\partial \bar{\gamma}}{ \partial p} \frac{\text{d} p}{\text{d} t}
\label{chap3_eq:gamma_partial_t}
\end{align}
The energy balance for each wall zone~\eqref{chap3_eq:WTPSC} reads,
\begin{align}
A_w\rho_wc_{p_w}\left(l_i\frac{\text{d}T_{w_i}}{\text{d}t}+\left(T_{w,B_{i,i-1}}-T_{w_i}\right)\frac{\text{d}z_{a,i}}{\text{d}t}+\left(T_{w_i}-T_{w,B_{i,i+1}}\right)\frac{\text{d}z_{b,i}}{\text{d}t}\right)\nonumber\\=\dot{Q}_{exh_i}-b_{WF}\alpha_{WF,i}l_i\left(T_{w_i}-\bar{T}_i\right) - \alpha_{amb}p_{evap}l_i\left(T_{w_i}-T_{amb}\right),
\label{chap3_eq:WTPSC}
\end{align}
where $A_w$, $\rho_w$ and $c_{p_w}$ are the wall cross-sectional area, density and heat capacity.
$T_{w_i}$ is the temperature of the respective wall zone and $T_{w,B_{i,i-1}}$ and $T_{w,B_{i,i+1}}$ are the wall temperatures at the left- and right-hand boundary of the zone, which are calculated using a length-weighted average, as suggested in \cite{Zhang2006}.
$\dot{Q}_{exh_i}$ is the amount of heat transfered from the exhaust gas to the wall, $b_{WF}$ is the WF channel width and $\alpha_i$ the heat transfer coefficient for the WF in the respective zone.
We introduce a term accounting for heat loss from the exchanger wall to the environment, in which $\alpha_{amb}$ is the heat transfer coefficient, $p_{evap}$ the HX perimeter and $T_{amb}$ the ambient temperature.
\\
\\
By analytical integration of the quasi-stationary energy balance on the exhaust side from interface $i+1$ to interface $i$ assuming static one dimensional flow \cite{McKinley2008}, the temperature at the end of one element can be calculated as in \eqref{chap3_eq:ExhaustTemp-i} and the heat transfered to the wall as in \eqref{chap3_eq:Qexh}.
Both the exhaust heat capacity $c_{p,exh_i}$ and the heat transfer coefficient $\alpha_{exh_i}$ are assumed constant over one element.
$\dot{m}_{exh}$, $T_{exh}$ and $b_{exh}$ are mass flow, temperature and width of the exhaust channel.
\begin{align}
T_{exh_i}&=T_{w_i}+\left(T_{exh_{i+1}}-T_{w_i}\right)\exp\left(-\frac{\alpha_{exh_i}b_{exh}}{\dot{m}_{exh}c_{p,exh_i}}l_i\right) ,~~~i \in [0,2]
\label{chap3_eq:ExhaustTemp-i}
\\
\dot{Q}_{exh_i}&=\dot{m}_{exh}c_{p,exh_i}\left(T_{exh_{i+1}}-T_{exh_i}\right)
\label{chap3_eq:Qexh}
\end{align}
\subsection{Pump and turbine models}
We model the pump assuming a fixed isentropic and mechanical efficiency ($\eta_{is,pump}$, $\eta_{mech,pump}$), according to \eqref{chap3_eq:Pump_Model}.
Within the model, we set both efficiencies to 0.9.
\begin{align}
P_{pump} = \frac{1}{\eta_{mech,pump}} \cdot \dot{m}_{WF} \cdot \frac{h_{out,is} - h_{in}}{\eta_{is,pump}}
\label{chap3_eq:Pump_Model}
\end{align}
For the turbine, we use \eqref{chap3_eq:Turb_Model} to calculate the power output $P_{turb}$.
\begin{align}
P_{turb} = \eta_{mech,turb} \cdot \dot{m}_{WF} \cdot \eta_{is,turb} \cdot \left( h_{in} - h_{out,is} \right) 
\label{chap3_eq:Turb_Model}
\end{align}
The isentropic efficiency is a function of pressure ratio between high and low pressure and turbine speed.
For this, we choose a polynomial function of third order with respect to pressure ratio and fifth order with respect to turbine speed.
The mechanical efficiency, in contrast, is a function of turbine speed $n$ (second order polynomial) and torque $M$ (fifth order polynomial).








\end{document}